\documentclass[onefignum,onetabnum]{siamart190516}
\usepackage{epstopdf}
\usepackage{algorithmic}
\usepackage{nicefrac}

\ifpdf%
    \DeclareGraphicsExtensions{.eps,.pdf,.png,.jpg}
\else
    \DeclareGraphicsExtensions{.eps}
\fi

\usepackage{amsopn}

\usepackage{booktabs}

\title{Reducing operator complexity in Algebraic Multigrid with Machine Learning Approaches}

\headers{ML approach for non-Galerkin coarse-grid operator}{R. Huang, K. Chang, H. He, R. Li and Y. Xi}

\author{Ru Huang\thanks{Department of Mathematics, Emory University, Atlanta, GA 30322 (\texttt{\email{\{ru.huang,kai.chang,yxi26\}@emory.edu}}). The work is supported by NSF awards  OAC 2003720 and  DMS 2208412.}  
 \and Kai Chang\footnotemark[1] 
\and Huan He\thanks{Department of Biomedical Informatics, Harvard University, Boston, MA 02130 (\texttt{\email{huan\_he@hms.harvard.edu}})}
    \and Rui Peng Li\thanks{Center  for Applied  Scientific Computing,
    Lawrence  Livermore National  Laboratory,  P. O.  Box 808,  L-561,
    Livermore,   CA  94551   {(\texttt{\email{li50@llnl.gov}})}.  This   work  was
    performed under the  auspices of the U.S. Department  of Energy by
    Lawrence    Livermore   National    Laboratory   under    Contract
    DE-AC52-07NA27344 and was supported by the LLNL-LDRD program under Project No. 23-FS-031.}
    \and Yuanzhe Xi\footnotemark[1]}

\usepackage{bm}
\usepackage{rotating}
\usepackage{color}
\usepackage{amsfonts}
\usepackage{listings}
\usepackage{float}
\usepackage{algorithm}
\usepackage{tikz}
\usepackage{mathrsfs}

\usepackage{enumitem}
\newlist{todolist}{itemize}{2}
\setlist[todolist]{label=$\square$}
\usepackage{pifont}

\def\Re{\mathbb{R}}

\def\Ee{\mathbb{E}}

\def\inv{^{-1}}
\def\invt{^{-\mathsf{T}}}
\def\trans{^{\mathsf{T}}}
\providecommand{\norm}[1]{\lVert#1\rVert}

\mathchardef\mhyphen="2D

\definecolor{codegreen}{rgb}{0,0.6,0}
\definecolor{codegray}{rgb}{0.5,0.5,0.5}
\definecolor{codepurple}{RGB}{50,157,168}
\definecolor{mycolor}{RGB}{51,177,255}
\definecolor{backcolour}{gray}{0.95}%{rgb}{0.95,0.95,0.92}

\lstdefinestyle{mystyle}{
  backgroundcolor=\color{backcolour},   commentstyle=\color{codegreen},
  keywordstyle=\color{mycolor},%magenta},
  numberstyle=\tiny\color{codegray},
  stringstyle=\color{codepurple},
  basicstyle=\ttfamily\footnotesize,
  breakatwhitespace=false,         
  breaklines=true,                 
  captionpos=b,                    
  keepspaces=true,                 
  numbers=left,                    
  numbersep=5pt,                  
  showspaces=false,                
  showstringspaces=false,
  showtabs=false,                  
  tabsize=2
}

\usepackage{multirow}

\newsiamremark{remark}{Remark}
\newsiamremark{hypothesis}{Hypothesis}
\crefname{hypothesis}{Hypothesis}{Hypotheses}
\newsiamthm{claim}{Claim}

\begin{document}

\maketitle

% \begin{center}
% In collaboration with:
%   {\TheCollaborators}
% \end{center}
% \vspace{1cm}

% \begin{todolist}
%     \item[\done] add V-cycle

%     \item make a new illustration figure

%     \item Add a new algorithm: sparsified V-cycle

%     \item[\done] modify abstract 

%     \item[\done] Modify introduction

%     % \item Question: how to justify the connection between minimizing the current loss function and the quantity $\phi$

%     \item 2-level and 3-level results
% \end{todolist}

\begin{abstract}
    % Algebraic multigrid (AMG) is an extensively studied scalable solver for solving symmetric positive definite linear systems. Although it has been shown, both practically and theoretically, to work very well for many problems, it certainly is not perfect. One associated issue lies in the gradually increasing density in the coarse grid matrices, meaning that as the level gets deeper, the matrix at each level has more non-zero entries. This would lead to increasing costs in communication, also known as data moving between processors, which can limit the parallel scalability of the method.
We propose a data-driven and machine-learning-based approach to compute 
non-Galerkin coarse-grid operators in algebraic multigrid (AMG) 
methods, addressing the well-known issue of increasing operator 
complexity. Guided by the AMG theory on spectrally equivalent 
coarse-grid operators, we have developed novel ML algorithms that 
utilize neural networks (NNs) combined with smooth test vectors from 
multigrid eigenvalue problems. The proposed method demonstrates promise 
in reducing the complexity of coarse-grid operators while maintaining 
overall AMG convergence for solving parametric partial differential 
equation (PDE) problems.  
Numerical experiments on 
anisotropic rotated Laplacian and linear elasticity problems 
are provided
to 
showcase the performance and compare with existing methods for 
computing non-Galerkin coarse-grid operators.
\end{abstract}

\begin{keywords}
  machine learning, multigrid methods,
operator complexity, neural networks
\end{keywords}
\begin{AMS}
65M55, 65F08, 65F10, 15A60
\end{AMS}

\section{Introduction}\label{sec-intro}

\begin{figure}[h!] \label{sparsity-pattern}
    \centerline{\includegraphics[width=0.5\columnwidth]{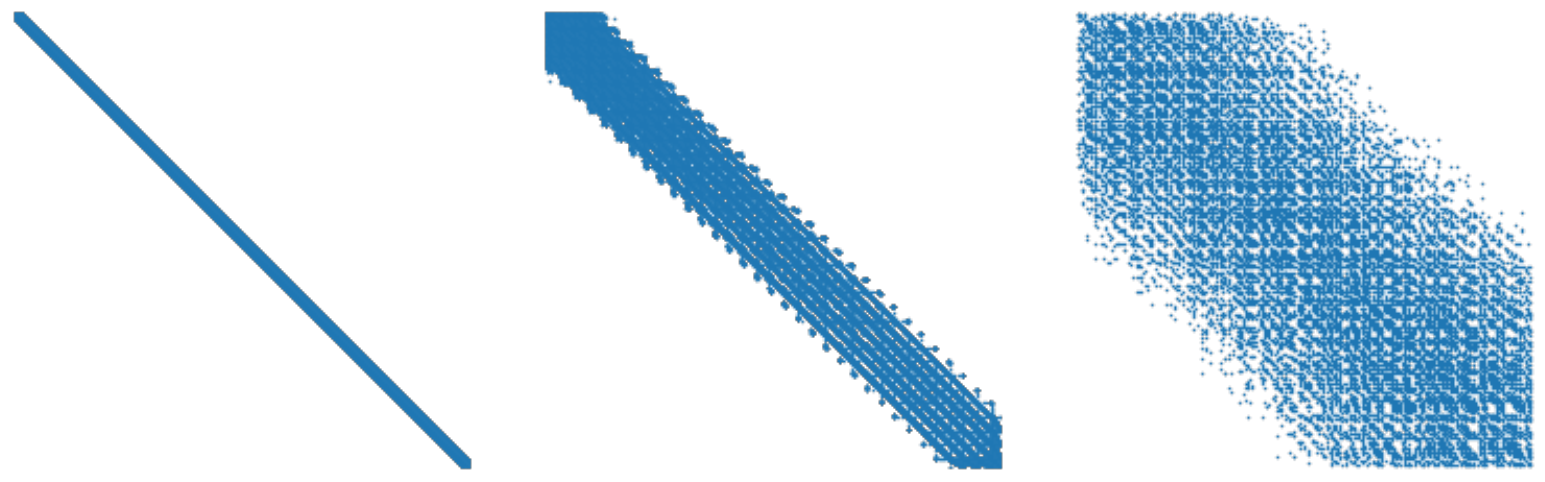}}
    \vspace{1em}
    \centering
    \begin{tabular}{c|c|c|c}
       level  & N & NNZ & RNZ\\\hline
        0&1,000,000 & 6,940,000 & 7\\
        1&500,000 &8,379,408 &17 \\
        2&128,630 &5,814,096 &45 \\
        3&23,023 &1,727,541 & 75\\
        4 & 11,688 & 1,371,218 & 117
        % 4&65025 & & \\
    \end{tabular}
    \caption{The sparsity patterns of $A^{(0)}$, $A^{(3)}$, and $A^{(5)}$ in the AMG hierarchy for solving the 3-D Poisson's equation (top). The size of the operator matrix (N), the number of nonzeros (NNZ) and the average number of nonzeros per row (RNZ) across the AMG levels (bottom).}
\end{figure}

\begin{figure}[h!]
    \centering
    \includegraphics[width=0.5\columnwidth]{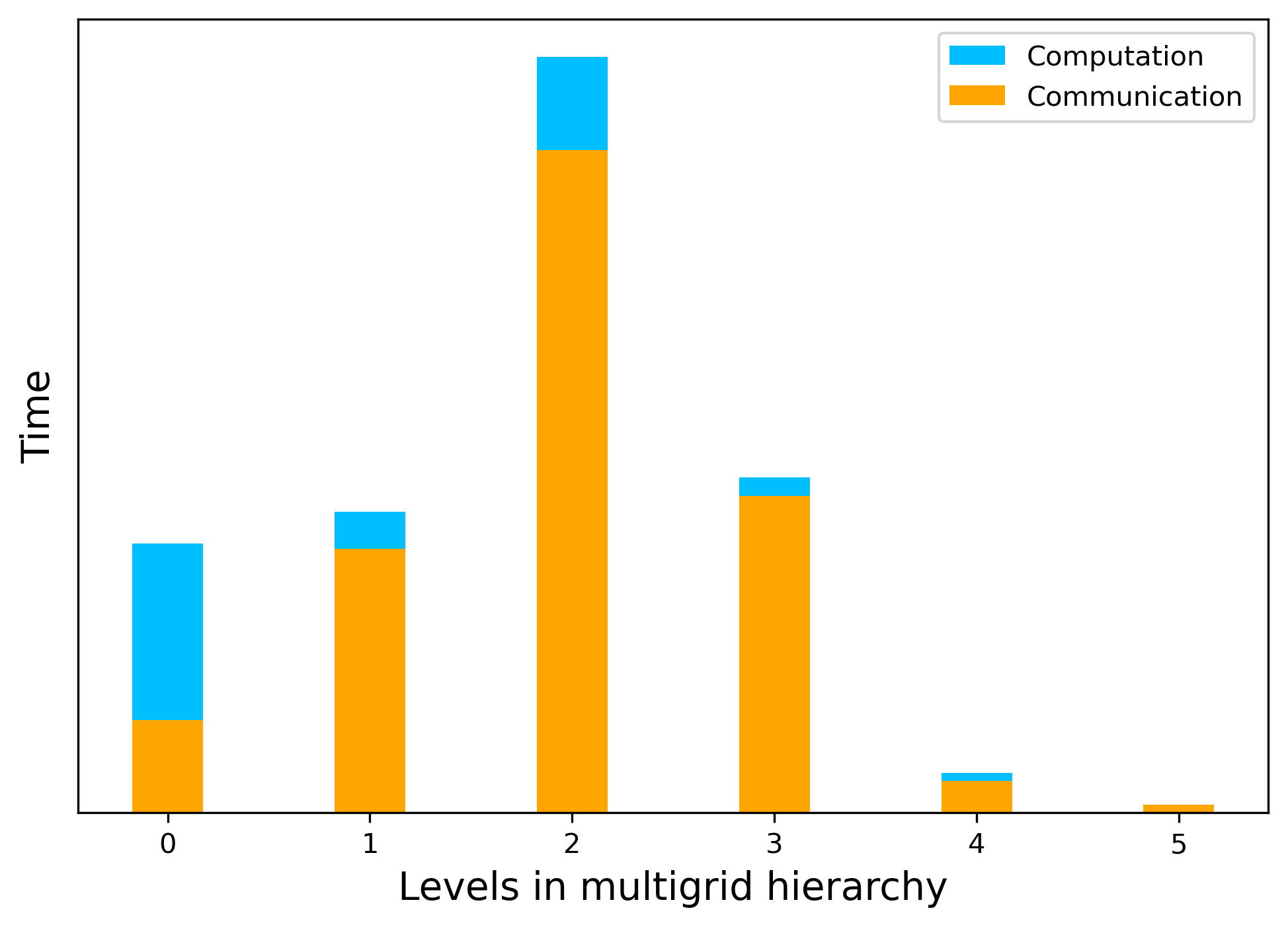}
    \caption{The cost of computation and communication in the first 6 levels of parallel AMG methods for solving the 3-D Poisson's equation. Image source: \cite{bienz2020reducing}}
    \label{communication}
\end{figure}

Algebraic Multigrid (AMG) methods are one of the most efficient and scalable iterative methods for solving linear systems of equations
\begin{equation}
\label{eq:1st-level}
    Au = f 
\end{equation}
where the coefficient matrix $A \in \Re^{N\times N}$ is  
sparse and large,
and $u\in \Re^{N}$ and $f\in \Re^{N}$ are the solution and
right-hand-side vectors respectively. 
For the systems that arise from 
elliptic-type partial differential equations (PDEs),  
AMG methods often exhibit optimal linear computational complexities.
Nevertheless, there is ongoing research focused on further improving the efficiency and scalability of AMG methods, in particular
for large-scale and challenging problems.
By and large,
the overall efficiency of iterative methods
is determined by not only 
the convergence rate of the iterations
but also the arithmetic complexity per iteration and the corresponding throughput on the underlying computing
platform. In this work, we address a common issue in AMG methods that is
the growth of the coarse-grid operator complexity in the hierarchy.
This operator is typically
computed as the Galerkin product from the operators in the 
fine level.
Assuming $A$ is symmetric positive definite (SPD),
the Galerkin operator is optimal  in the sense that
it yields an orthogonal
projector as the coarse-grid correction that guarantees to reduce
the $A$-norm of the error.
However, on the other hand, this operator can lead to the issue of increasing operator sparsity, particularly at deeper levels of the AMG hierarchy.
This can impair the overall performance of AMG by
introducing challenges in terms of computational efficiency, 
memory requirements,
and the communication cost in distributed computing environments.
Moreover, the increasing operator complexity can also affect the effectiveness and robustness of other AMG components 
such as the coarsening and interpolation algorithms.
To demonstrate this problem, we consider classical
AMG methods for solving the 3-D Poisson's equation discretized 
on a $100\times 100\times 100$ grid with a 7-point stencil.
The sparsity patterns of the operator matrix $A^{(l)}$ at the levels $l=0,3,5$
are shown in  \cref{sparsity-pattern}. From these patterns, it is evident that
the matrix bandwidth  increases as the level goes deeper,
as well as the stencil size (i.e., the average number of nonzeros per row).
The increased sparsity often leads to not only a growth in computational cost but also an increase in data movement, which corresponds to the communication expense in parallel solvers.
\cref{communication} shows the time spent
in the computation and communication in the first 6 levels of the 
AMG hierarchies for solving this 3-D Poisson's problem.
As depicted, there is a steep increase in the computational cost at level 2, 
coinciding with the level where the communication cost reaches its maximum.

One approach to reducing the coarse-grid operator complexity is to
``sparsify'' the
Galerkin operator after it is computed, 
i.e., removing some nonzeros outside a given sparsity
pattern. The obtained sparsified operator is often called a ``non-Galerkin''
coarse-grid operator. The methods developed in \cite{wienands2009collocation,treister2012algebraic} leverage algebraically smooth basis vectors and  the
approximations to the fine grid operator 
to explicitly control the coarse grid sparsity pattern.
The algorithms introduced in \cite{falgout2014non,treister2015non} first determine the patterns of the sparsified operator based on 
heuristics on the path of edges in the corresponding graph
and then compute the 
numerical values to ensure the  spectral equivalence to the Galerkin operator
for certain types of PDEs. 
 Improving the parallel efficiency of AMG by reducing the communication cost with the non-Galerkin operator was discussed in \cite{bienz2016reducing}. These existing algorithms for computing non-Galerkin operators  are usually
based on heuristics on the associated graph and the characteristic
of the underlying PDE problem, such as the information of the near kernels of $A$.
Therefore, they are problem-dependent, and often times it can be 
difficult to devise such heuristics that are suitable for a broader class of problems.

Recently, there has been a line of work in the literature to leverage data-driven and ML-based methods to improve the robustness of AMG. 
In particular, \cite{luz2020learning, katrutsa2020black, greenfeld2019learning} deal with learning better prolongation operators.  
Techniques of deep reinforcement learning (DRL) 
are exploited in \cite{taghibakhshi2021optimization}
to better tackle the problem of AMG coarsening
combined with the diagonal dominance ratio of the F-F block. Both the works in \cite{huang2021learning} and \cite{kuznichov2022learning} focus on the problem of designing better smoothers. In \cite{huang2021learning}, smoothers are directly parameterized by multi-layer 
convolution neural networks
(CNNs) while \cite{kuznichov2022learning} optimizes the weights in the weighted Jacobi smoothers. In this paper, we follow this line of research and propose a data-driven and ML-based method  
for non-Galerkin operators. 
The innovations and features of the proposed method are summarized as follows:
1) Introduction of a multi-level algorithm based on ML methods to sparsify all coarse-grid operators in the AMG hierarchy;
2) Successful reduction of operator density while preserving the convergence behavior of the employed AMG method;
3) Applicability of the proposed NN model to a class of parametric PDEs with parameters following specific probability distributions;
4) Ability to train the sparsified coarse-grid operator on each level in parallel once the training data is prepared;
5) Flexibility for the user to choose the average number of non-zero entries per row in the coarse-grid operators, with a minimum threshold requirement.
To the best of our knowledge, our proposed work is the first to utilize ML models for controlling sparsity within the hierarchy of AMG levels.

The rest of the paper is organized as follows. We first briefly review the preliminaries of AMG methods and the non-Galerkin algorithms 
in \cref{sec-mg}. We elaborate on our proposed sparsification algorithm in \cref{sec-spdl}. Numerical experiments and results are presented in \cref{sec-numerical-exp}. Finally, we conclude in \cref{sec-conclusion}.

% In this work, we tackle the problem of increasing complexities of coarse-grid operators, computed as Galerkin products,
% along the levels of AMG hierarchies. This problem has adverse effects of the overall performance of AMG as it can increase the computational cost of applying the operators in deeper levels, can impair the effectiveness and robustness of other AMG components such as coarsening algorithms and interpolation algorithms, and can make the communications more expensive in distributed computing environments, to name a few.
% Traditional ``Non-Galerkin'' approaches to reduce the operator complexity employs
% sparsification of the Galerkin product [More to write].
%In this work, we will limit our discussion to the multigrid V-cycle as it is the most natural way to study MG. 
%We first give a short review on the V-cycle scheme along with introducing our notations.

% \input{tex/1.1-vcycle.tex}
% \input{tex/1.2-related.tex}
\section{AMG preliminaries and coarse-grid operators}\label{sec-mg}

% Multigrid methods exploit a hierarchy of grids to speed up the convergence of relaxation methods. Other than a relaxation method, the major components of a basic multigrid method consist of a restriction operator $\mxR$ and a prolongation operator $\mxP$. The restriction operator is used for transferring information from levels with finer grids (a.k.a. fine levels) to levels with coarser grids (a.k.a. coarse levels) while the prolongation operator is used for transferring information from coarse levels to fine levels. We first discuss relaxation methods, introduce algebraic multigrid (AMG) along with restriction and prolongation operators, and demonstrate the issue of increasing communication costs in multigrid methods.

% Our remedy for the aforementioned issue of increasingly dense coarse-level operators is based on the notion of spectrally equivalent stencils. 
% Operators generated by spectrally equivalent stencils have similar spectrum which determines the convergence behavior of MG. Previous results have provided ways of quantifying the spectral equivalence between stencils. 
%Before we introduce our method and the issue that we want to address, we will 
%first provide a brief overview of Multigrid and its theory.
In this section, we give a brief introduction to AMG methods 
and the Galerkin coarse-grid operators. 
The AMG method is a multilevel method that utilizes a hierarchy of grids, 
consisting of fine and coarse levels, 
and constructs coarse-level systems at different scales that can capture 
the essential information 
of the fine-level system while reducing the problem size. 
AMG algorithms employ techniques such as coarsening, relaxation, 
restriction and interpolation
to transfer information between the grid levels 
to accelerate the solution process.
 \cref{alg:v} presents the most commonly used AMG V-cycle scheme.
It uses $\nu$ steps of pre- and post-smoothing, where $M$
and $M\trans$ are the smoothing operators.
Matrices $R$ and $P$ are the restriction and prolongation operators, respectively.
The coarse-grid operator  
is computed in Step \ref{alg:VGalerkin} via the Galerkin product,
$A_g=RAP$.
The aim of the smoothings is to quickly annihilate the high-frequency 
errors via simple
iterative methods such as relaxation, whereas the low-frequency errors
are targeted by the Coarse-Grid Correction (CGC) operator, $I-P(RAP)\inv R A$.
When $R=P\trans$, the CGC operator is
$A$-orthogonal with the Galerkin operator $A_g$. 

\begin{algorithm}[h!]
    \caption{Multigrid V-Cycle for solving $A^{(l)}u^{(l)}=f^{(l)}$ at level $l$}
    \label{alg:v}
    \begin{algorithmic}[1]
        %\REQUIRE $l$, $L$, $\mxA_{\beta}$, $\bff$, $\{\mxP^{(k)}\}_{k=1}^{L-1}$, $\{\mxR^{(k)}\}_{k=1}^{L-1}$, $\mxM$, $\nu$    
        % \hspace{1px}
        %\ENSURE $\bfu^{(l)}$ such that $\mxA_{\beta}\bfu \approx \bff$ 
        % \hspace{1px}
        % \hrule
        % \hspace{1px}
        %\STATE Initialize $\bfu^{(l)} = \bfzero$ 
        \STATE Pre-smoothing: $u^{(l)} := (I-(M^{(l)})^{-1} A^{(l)}) u^{(l)} + (M^{(l)})^{-1}f^{(l)}$ for $\nu$ steps \label{alg:Vpresmooth}
        \STATE Compute residual $r^{(l)} = f^{(l)} - A^{(l)}u^{(l)}$ and the restriction
        $r^{(l+1)} = R^{(l)}r^{(l)}$
        % \State Let $\mxA_{g,\beta}^{(0)} = \mxA$ and $\bfr_g^{(0)} = \bfr$}
        \STATE Compute Galerkin operator $A_{g} = R^{(l)}A^{(l)}P^{(l)}$ and let $A^{(l+1)}=A_g$
        \label{alg:VGalerkin}
        \IF{$l = L-1$}
            \STATE Solve $A^{(l+1)}u^{(l+1)} = r^{(l+1)}$ with an arbitrary method
            % \State \Return $\bfe_g^{(l+1)}$
        \ELSE
            \STATE Let $u^{(l+1)}=0$ and $f^{(l+1)}=r^{(l+1)}$. Go to Step \ref{alg:Vpresmooth}
            with $l:=l+1$.
        \ENDIF
        \STATE Prolongate and correct: $u^{(l)} := u^{(l)}+P^{(l)}u^{(l+1)}$ 
        % \For{$i = 1:\nu$}
        % \State $\bfu \gets (\mxI-\mxM^{-1} \mxA) \bfu + \mxM^{-1} \bff $
        % \EndFor
        \STATE Post-smoothings: $u^{(l)} := (I-(M^{(l)})\invt A^{(l)}) u^{(l)} + (M^{(l)})\invt f^{(l)} $ for $\nu$ steps 
        %\RETURN $\bfu^{(l)}$ 
    \end{algorithmic}
\end{algorithm}

\subsection{Non-Galerkin  operators}
Naive approaches, such as indiscriminately removing nonzero entries in the Galerkin operators based on the magnitude often result in slow convergence of the overall AMG method (see the example provided in Section 3 of \cite{falgout2014non}). To address the aforementioned challenges  arising from the increased operator complexity in the Galerkin operator $A_g$, alternative operators, denoted by $A_c$, that are not only
sparser than $A_g$ but also spectrally equivalent
have been studied and used in lieu of  the Galerkin operator.

\begin{definition}
SPD matrices  $A_g$ and $A_c$ are spectrally equivalent if
\begin{equation} \label{eq:speq}
    0 < \alpha \le \lambda(A_g\inv A_c) \le \beta,  
\end{equation}
with $\alpha$ and $\beta$ both close to 1.
\end{definition}
The convergence rate of AMG can be analyzed through the spectral radius of the error propagation matrix. For example, the two-grid error propagation matrix corresponding to the V-cycle in \cref{alg:v} reads 
\begin{equation} \label{eq:TGE_G}
E_{g} = (I-M\invt A)^{\nu}(I-PA_{g}^{-1} R A) (I-M^{-1} A)^{\nu}.
\end{equation}
With the replacement of $A_g$ by $A_c$, it becomes
\begin{equation} \label{eq:TGE_C}
E_{c} = (I-M\invt A)^{\nu}(I-PA_{c}^{-1} R A) (I-M^{-1} A)^{\nu}.
\end{equation}
The spectrum property of $E_{c}$ is analyzed in the following theorem.

\begin{theorem}[\cite{falgout2014non}]
\label{thm:theta}
Denoting by $B_g$ and $B_c$ respectively the corresponding preconditioning matrices 
induced by $E_g = I-B_g^{-1}A$
and $E_c=I-B_c^{-1}A$, and assuming $A_{c}$ and $A_{g}$ are both SPD and
\begin{equation}
\label{eq:theta}
\theta=\|I-A_{c}A_{g}^{-1}\|_{2} = \|I-A_{g}^{-1}A_{c}\|_{2} < 1,
\end{equation}
for the preconditioned matrix, we have
\begin{equation}
\kappa(B^{-1}_cA)\leq \frac{1+\theta}{1-\theta}\kappa(B^{-1}_gA),    
\end{equation}
and moreover
\begin{equation}
\label{eq:ec-bound}
    \rho(E_c)\leq \max\left(\frac{\lambda_{\max}(B^{-1}_gA)}{1-\theta}-1,1-\frac{\lambda_{\min}(B_g^{-1}A)}{1+\theta}\right),
\end{equation}
where $\lambda_{\mathrm{max}}(\cdot)$ and $\lambda_{\mathrm{min}}(\cdot)$ are the largest and smallest eigenvalues respectively, $\kappa(\cdot)$ denotes the condition number and $\rho(\cdot)$ denotes the spectrum radius.
\end{theorem}
The quantity $\theta$ measures the degree of
spectral equivalence between the operators $A_g$ and $A_c$, i.e., 
only when $\theta$ is small, these operators are 
spectrally equivalent.
Clearly, the condition number of the preconditioned matrix and
the two-grid convergence with respect to $E_c$
deteriorate  as $\theta$ increases.
With fixed $B_g\inv A$, 
we can establish a criterion for the convergence of $E_c$ 
with respect to $\theta$, as shown in the next result.

\begin{corollary}\label{cor:theta}
Suppose $\theta<1-\lambda_{\max}(B_g^{-1}A)/2$, the two-grid method \cref{eq:TGE_C} converges. 
% Moreover, if TG with $\mxA_{g}$ converges and $\phi<1/2$, then the method with $\mxA_{c}$ is guaranteed to converge.
\end{corollary}

\begin{proof}
    Note that 
    $\theta < 1 - \lambda_{\max}(B_g^{-1}A)/{2}$
    implies that ${\lambda_{\max}(B_g^{-1}A)}/{(1-\theta)} - 1 < 1$. 
    Since $\lambda(B_g\inv A) > 0$, 
    $1-\lambda_{\min}(B_g^{-1}A)/{(1+\theta)} < 1$.
Therefore, from \eqref{eq:ec-bound}, it follows that $\rho(E_c) < 1$.
\end{proof}

\subsection{Spectrally equivalent stencils}
In this paper, we focus on structured matrices that can be represented by stencils and grids. These structured matrices exhibit a unique property wherein two spectrally equivalent stencils can determine two sequences of spectrally equivalent matrices with increasing sizes and thus ensures the convergence of $E_c$ with increasing matrix sizes.
\begin{definition}[\cite{axelsson1982multigrid, bolten2007structured}]
\label{def:spectral-equivalence}
{Let $\{A_{j}\}$ and $\{B_{j}\}$ be two sequences of (positive definite) matrices with increasing size $N_{j}$, where $A_j$ and $ B_j \in \mathbb{R}^{N_j\times N_{j}}$. 
%If the eigenvalues of $B_{j}^{-1}A_{j}$, $\lambda(B_{j}^{-1}A_{j})$, satisfy
%$$0<\alpha<\lambda(B_{j}^{-1}A_{j})\leq \beta<\infty$$
If $A_j$ and $B_j$ are spectrally equivalent as defined in \cref{eq:speq}
for all $j$ with $\alpha$ and $\beta$ that are independent of $N_j$, 
then the sequences $\{A_{j}\}$ and $\{B_{j}\}$ are called 
spectrally equivalent sequences of matrices.}
\end{definition}

%Intuitively, this definition formalizes the similarity between the spectrum of two sequences of matrices. This provides a natural way of defining spectrally equivalent stencils. Recall that a stencil can represent a matrix of any size, as long as the pattern remains the same. Therefore, for two stencils to be spectrally equivalent, it is reasonable to require the spectrum of any two matrices generated by the respective stencils to be similar, as long as the sizes of the two matrices are the same. This is summarized in the following definition.
The above definition yields the definition of spectrally equivalent stencils. 
\begin{definition}[\cite{bolten2007structured}]
Suppose the sequences of matrices $\{A_{j}\}$ and $\{B_{j}\}$ 
are constructed with the stencils $\mathcal{A}$ and
$\mathcal{B}$ respectively, where $A_j$ and $B_j$ have the same size for all $j$. 
We call $\mathcal{A}$ and $\mathcal{B}$ are spectrally equivalent 
if $\{A_{j}\}$ and $\{B_{j}\}$ are spectrally equivalent sequences of matrices.
\end{definition}

%With the definition, several questions could be of immediate interest: 1) does there exist spectrally equivalent stencils? 2) for a fixed stencil, does there exist a stencil spectrally equivalent to that particular one? 3) if the answer to the last question is positive, how to find such a stencil computationally? and 4) is it possible for the spectrally equivalent stencil to have fewer non-zero entries so that the multigrid method can be benefited?

%We first consider a special class of stencils, namely the circulant stencils (which will be defined shortly). In fact, it has been shown that there exists a closed-form 5-point stencil that is spectrally equivalent to a 9-point circulant stencil \cite{bolten2007structured}. We summarize this in the next theorem.

At the end of this section, we provide an example of spectrally equivalent stencils. Consider the following 9-point stencil that was used in the study of the
AMG method for circulant matrices \cite{bolten2007structured}:
\begin{equation} \label{eq:9ptcirculant}
\left[
\begin{array}{ccc}
c&b&c\\
a&-2(a+b)-4c&a\\
c&b&c\\
\end{array}
\right].
\end{equation}
It was proved that the associated 5-point stencil
\begin{equation} \label{eq:5ptcirculant}
\left[
\begin{array}{ccc}
&b+2c&\\
a+2c&-2(a+b)-8c&a+2c\\
&b+2c&\\
\end{array}
\right].
\end{equation}
is spectrally equivalent to \cref{eq:9ptcirculant}.
Results for the 7-point stencil in 3-D that is spectrally equivalent to a 
27-point stencil can be also found in \cite{bolten2007structured}.

\subsection{Numerical heuristics for spectral equivalence}
\label{sec:num_spec_eq}
Directly optimizing \cref{eq:theta} to find a spectrally equivalent
$A_c$ appears to be challenging. Instead, a more viable approach is to utilize test vectors that correspond to the low-frequency modes of $A_g$ (see, e.g., \cite{falgout2014non,wienands2009collocation}).
These low-frequency modes represent the algebraically 
smooth modes at a coarse level, which are important for the interpolation to 
transfer to the fine level within the AMG hierarchy.
From
the perturbed error propagation operator
\cref{eq:TGE_C}, it follows that
after the pre-smoothing steps, the remaining error, denoted by $e$,
that is algebraically smooth in terms of
$A$ (i.e., $Ae \approx 0$)
needs to be efficiently annihilated by the coarse-grid.
By the construction of the interpolation operator, 
the smoothed error
is in the range of 
$P$, meaning that, 
$e = Pe_c$ with some coarse-grid error $e_c$.
Furthermore, $e_c$ is smooth with respect to
$A_g$, since $A_g e_c =RAPe_c = RAe$ is small.
Therefore, for an effective CGC with non-Galerkin 
coarse-grid operator $A_c$, 
it is essential for
%\begin{align}
$(I-PA_c\inv RA)e = (I-PA_c\inv RA)Pe_c =P(I-A_c\inv A_g)e_c$
%\notag
%\end{align}
to be small, which implies $A_g e_c \approx A_c e_c$.
That is to enforce the accuracy of the application of $A_c$ on low-frequencies vectors with respect to that of $A_g$.

% This is because, suppose that $A_g$ admits eigendecomposition 
%$A_g=U\Lambda U\trans$ and let $E= A_g - A_c$, 
%for any $x$ such that $\norm{x}_2=1$, we can write
%\begin{align}
%    \norm{(I-A_c A_g\inv) x}_2 &= \norm{E A_g\inv U y}_2
%    = \norm{E U\Lambda\inv y}_2 = \norm{\textstyle\sum_i ({y_i}/{\lambda_i})Eu_i}_2 \notag \\
%    &\le \textstyle\sum_i \vert y_i/\lambda_i\vert \norm{Eu_i}_2 
%    \le \textstyle\sum_i \norm{Eu_i}_2/\lambda_i, \notag
%\end{align}
%where $\lambda_i$'s are the eigenvalues of $A_g$
%and $y_i$ denotes the $i$-th entry of $y=U\trans x$.
%Therefore, to have a small $\theta$, it is important to make $A_g u_i \approx A_c u_i$ 
%small for the eigenvectors $u_i$ that correspond to the small eigenvalues, i.e., 
%enforce the accuracy of the application of $A_c$ on low-frequencies modes with respect to that of $A_g$.

In this paper, we adopt the approach of multigrid eigensolver (MGE)
\cite{doi:10.1137/090752973} to compute
the smooth vectors in the AMG hierarchy. 
First, consider two-grid AMG methods. The Rayleigh quotient of  $Pe_c$ with respect to $A$
reads
\begin{equation} \label{eq:rayleigh}
   r(A, Pe_c)= \frac{(APe_c, Pe_c)}{(Pe_c, Pe_c)} = \frac{(P\trans A Pe_c, e_c)}{(P\trans Pe_c, e_c)} = \frac{(A_g e_c, e_c)}{(T e_c, e_c)} ,
\end{equation}
where $T=P\trans P$.
Therefore, the desired smooth modes that minimize \cref{eq:rayleigh} relate to
the eigenvectors that correspond to the
small eigenvalues of the generalized eigenvalue problem
\begin{equation}
    A_g u = \lambda T u , \quad T=P\trans P.
    \label{eq:smoothvec}
\end{equation}
For AMG methods with more than 2 levels, we can compute the smooth vectors 
at each coarse level by recursively applying
\cref{eq:smoothvec} at the previous fine level.

\section{An ML method for coarse-grid operators}
\label{sec-spdl}
%An approach to alleviate the problem of having rapidly increasing operator 
%complexities along the AMG multilevel hierarchy is
%to sparsify the coarse-grid operators after the Galerkin product~\eqref{eq:galerkin}.
%Therefore, this type of approaches are typically referred to as
%``non-Galerkin'' methods and the sparsified matrix is often called
%a non-Galerkin coarse-grid operator.
%Therefore, when designing a non-Galerkin method, 
%it is crucial to not deteriorate the spectral equivalence of
%non-Galerkin operator, at least not too much. 

%\subsection{The Algorithmic Pipeline}

\begin{figure}[h!]
\centering
\includegraphics[width=0.75\textwidth]{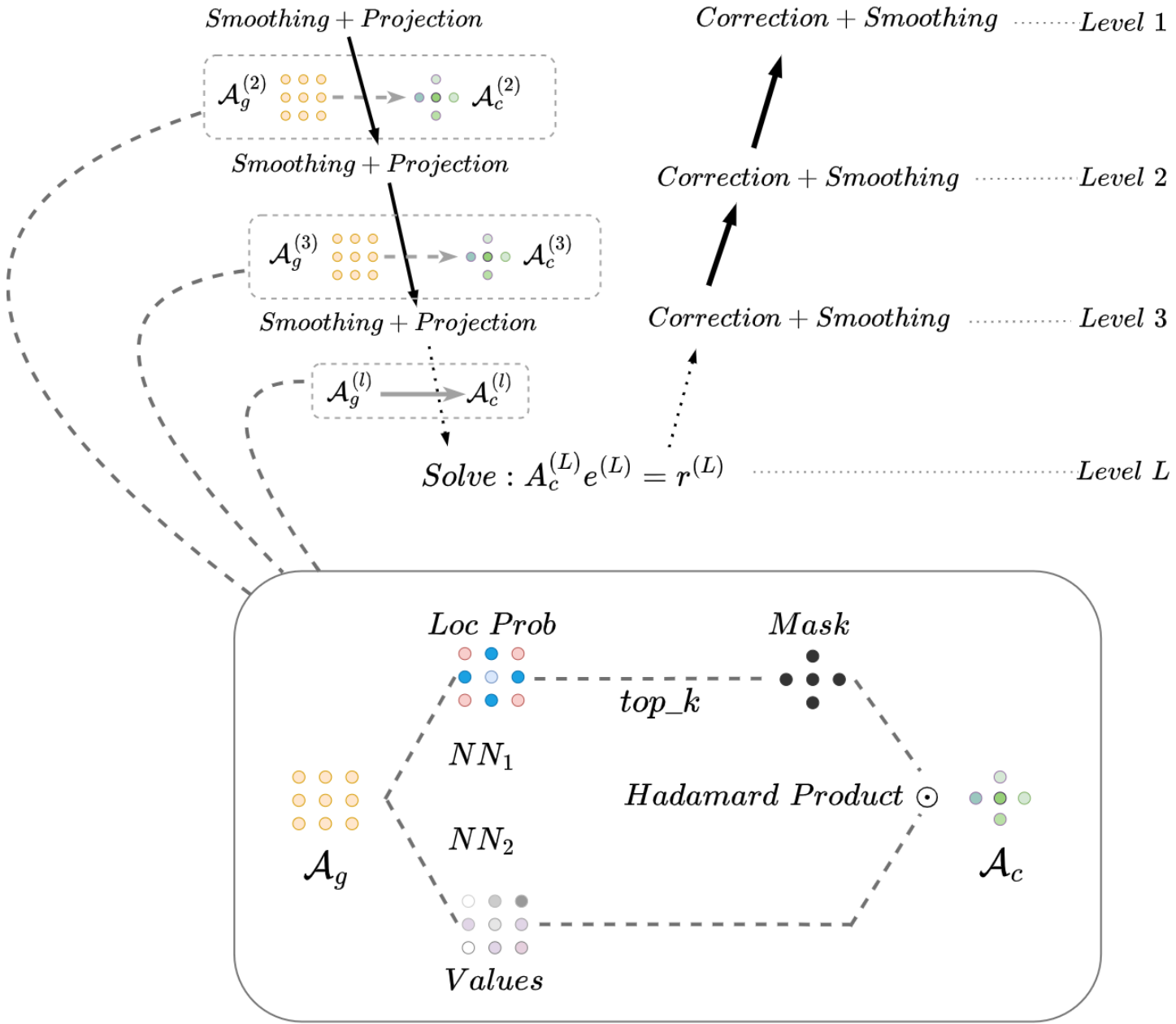}
\caption[Framework]{Illustration of the ML algorithm for
computing coarse-grid operators with NNs.}
\label{fig:algo-illustration}
\end{figure}

% \textbf{HOW TO PROCEED? Introduce algorithmic pipeline first or introduce the training first?}
\label{sec:ML}
We aim to utilize ML techniques to compute non-Galerkin operators 
in the AMG method for solving \cref{eq:1st-level},
where $A$ is a stencil-based coefficient matrix that
corresponds to PDE problems discretized on Cartesian grids.
%For the matrix stencil, denoted by $\beta \in \Re^p$, we assume
%$\beta \in \mathscr{B}$ that is associated with probability distribution 
%$p_{\beta}$. 
On a given AMG level $l>1$, with  stencil $\mathcal{A}_{g}^{(l)}$ 
associated with the Galerkin matrix $A_g^{(l)}$, 
we construct a sparser stencil $\mathcal{A}_{c}^{(l)}$  
in the following 3 steps: 
\begin{enumerate}[leftmargin=4.5em]
    \item[Step 1:] Select the pattern of $\mathcal{A}_{c}^{(l)}$, i.e., 
    the positions of the nonzero entries, where the corresponding entries 
    %in $\mathcal{A}_{g}^{(l)}$ 
    are 
    assumed to be nonzero;
    \item[Step 2:] Compute the numerical values of the nonzero entries;
    \item[Step 3:] Construct $\mathcal{A}_{c}^{(l)}$ by the point-wise multiplication
    of the pattern and values,
\end{enumerate}
which are explained in detail below 
and illustrated in \cref{fig:algo-illustration}.
\paragraph{Step 1} The NN in this step, denoted by $F_{\theta^{(l)}}$, is parametrized by
$\theta^{(l)}$. 
It computes the \emph{location probability} for each of the
stencil entries of $\mathcal{A}_{g}^{(l)}$.
We apply the NN $F_{\theta^{(l)}}$ to the vectorized stencil 
$v_g^{(l)}=\mbox{vec}(\mathcal{A}_{g}^{(l)})$, i.e.,
the vector reshaped from the stencil array,
followed by a \emph{softmax} layer.
Therefore, the output of the NN can be written as
\begin{equation}
\mathcal{P}^{(l)} =
\mbox{softmax}(F_{\theta^{(l)}}(v_{g}^{(l)})),
\end{equation}
which is then reshaped back to match the shape of $\mathcal{A}_{g}^{(l)}$.
Given that $0<\mathcal{P}^{(l)}<1$, each entry can be interpreted as the probability of the entry appearing (being nonzero) in the 
sparsified stencil $\mathcal{A}_{c}^{(l)}$. 
With that, we select the largest $k$ entries of $\mathcal{P}^{(l)}$,
\begin{equation}
\mathcal{I}^{(l)} = \left\{i \, \Bigr\vert \,\mathcal{P}^{(l)}(i) \,\, \text{is one of the largest $k$ entries of } \mathcal{P}^{(l)}\right\},
\end{equation}
where $\mathcal{I}^{(l)}$ denotes the set of the indices of those entries,
from which
we build a mask Boolean vector $\mathcal{M}^{(l)}$ defined as 
\begin{equation}
\mathcal{M}^{(l)}\left(i\right) = 
\begin{cases}
1, \quad & \text{if} \,\, i \in \mathcal{I}^{(l)} \\
0, \quad & \text{otherwise}
\end{cases} \ ,
\end{equation}
that determines the positions of the nonzeros in the non-Galerkin stencil. 
\paragraph{Step 2} The NN in this step, denoted by $G_{\psi^{(l)}}$, is parametrized by $\psi^{(l)}$ 
 which is applied to the same input as in Step 1. The output from NN of this step reads
\begin{equation}
\mathcal{V}^{(l)} = G_{\psi^{(l)}} (v_{g}^{(l)}),
\end{equation}
which determines the numerical values of the nonzero entries.
\paragraph{Step 3} The final non-Galerkin stencil is computed by the Hadamard (or element-wise) product  
\begin{equation}
\mathcal{A}_{c}^{(l)} = \mathcal{M}^{(l)} \odot \mathcal{V}^{(l)}.
\end{equation}
We summarize these steps in \cref{alg:sparsify}.
The AMG V-cycle using the sparsified coarse grid is outlined in 
\cref{alg:sparse-v}, which closely resembles \cref{alg:v}, whereas, 
instead of using the Galerkin operator for coarser levels, the non-Galerkin operator $A_c$ is constructed from the sparsified stencil generated by \cref{alg:sparsify}.

%\textcolor{red}{Huan: do we need a sentence to clarify the connection between these 
%algorithms}
\begin{remark}
A few remarks on \cref{alg:sparsify} and \cref{alg:sparse-v} follow. 
To begin with, the parameter $k$ of \cref{alg:sparsify} signifies the number of the nonzero entries in the sparsified stencil. This effectively gives us the ability to directly manipulate the complexity of the resulting non-Galerkin operator.
Secondly, in the NN implementations, we ensure that the shapes of $\mathcal{M}^{(l)}$ and $\mathcal{V}^{(l)}$ are identical. This enables the proper application of the Hadamard product.
Lastly, it is assumed that the NNs, $F_{\theta_l}$ and $G_{\psi_l}$, have undergone sufficient training. Therefore, Step \ref{alg:sparse-v-NNstep} of \cref{alg:sparse-v} involves merely the application of the trained NNs. 
\end{remark}

%%%%% Sparsify 
\begin{algorithm}[h!]
\caption{SparsifyStencil}
\label{alg:sparsify}
\begin{algorithmic}[1]
\REQUIRE $\mathcal{A}_{g}, F_{\theta},G_{\psi},k$
 \STATE Apply the NNs to compute $\mathcal{P}=F_{\theta}(\mathcal{A}_{g})$ and $\mathcal{V}=G_{\psi}(\mathcal{A}_{g})$
 \STATE $\mathcal{M}$ has the same shape as $\mathcal{V}$ and has a value of 1 at the entries corresponding to the $k$ largest values of $\mathcal{P}$, with 0 elsewhere.
 \RETURN  $\mathcal{A}_{c}=\mathcal{M}\odot \mathcal{V}$
 \end{algorithmic}
\end{algorithm}

%%%% Sparsified V-Cycle
\begin{algorithm}[h!]
    \caption{AMG V-Cycle with sparsified coarse-grid operator}
    \label{alg:sparse-v}
    \begin{algorithmic}[1]
%        \REQUIRE $l$, $L$, $\mxA$, $\bff$, $\mxM$, $\nu$, $k$, $\{\mxP^{(l)}\}_{l=1}^{L-1}$, $\{\mxR^{(l)}\}_{l=1}^{L-1}$, $\{f_{\theta_l}^{(l)}\}_{l=2}^{L}$, $\{g_{\psi_l}^{(l)}\}_{l=2}^{L}$  
        % \hspace{1px}    
%        \ENSURE $\bfu^{(l)}$ such that $\mxA\bfu^{(l)} \approx \bff$ 
        % \hspace{1px}
        % \hrule
        % \hspace{1px}
        %\STATE Initialize $\bfu^{(l)} = \bfzero$
        \STATE Pre-smoothing: $u^{(l)} := (I-(M^{(l)})^{-1} A^{(l)}) u^{(l)} + (M^{(l)})^{-1}f^{(l)}$ for $\nu$ steps \label{alg:Vpresmooth22}
        \STATE Compute the residual $r^{(l)} = f^{(l)} - A^{(l)}u^{(l)}$ and restriction
        $r^{(l+1)} = R^{(l)}r^{(l)}$
        % \State Let $\mxA_{g,\beta}^{(0)} = \mxA$ and $\bfr_g^{(0)} = \bfr$}
        \STATE Compute the Galerkin operator $A_{g}^{(l+1)} = R^{(l)}A^{(l)}P^{(l)}$
        \IF{$l = L-1$}
            \STATE Solve $A_{g}^{(l+1)} u^{(l+1)} = r^{(l+1)}$ with an arbitrary method
            % \State \Return $\bfe_g^{(l+1)}$
        \ELSE
            \STATE Apply \cref{alg:sparsify}: $\mathcal{A}_c^{(l+1)} = \mbox{SparsifyStencil}(\mathcal{A}_g^{(l+1)}, F_{\theta^{(l+1)}},G_{\psi^{(l+1)}}, k)$
            \label{alg:sparse-v-NNstep}
            \STATE  Let $A^{(l+1)}=A_c^{(l+1)}$, $u^{(l+1)}=0$ and $f^{(l+1)}=r^{(l+1)}$. Go to Step \ref{alg:Vpresmooth22}
            with $l:=l+1$
        \ENDIF
        \STATE Prolongate and correct: $u^{(l)} := u^{(l)}+P^{(l)}u^{(l+1)}$ 
        % \For{$i = 1:\nu$}
        % \State $\bfu \gets (\mxI-\mxM^{-1} \mxA) \bfu + \mxM^{-1} \bff $
        % \EndFor
        \STATE Post-smoothings: $u^{(l)} := (I-(M^{(l)})\invt A^{(l)}) u^{(l)} + (M^{(l)})\invt f^{(l)} $ for $\nu$ steps \end{algorithmic}
\end{algorithm}

\subsection{NN training algorithm}
In this section, we delve into the specifics of the training algorithm that enables 
NNs to generate stencils with a higher level of sparsity compared to the
Galerkin stencil, without impairing the overall convergence of the AMG method.
A key component of  \cref{alg:sparse-v} is line 
\ref{alg:sparse-v-NNstep}
where $F_{\theta}$ and $G_{\psi}$ are the pre-trained NNs. 
The loss function, another crucial component, is pivotal to the training procedure. 
Based on the discussion in \cref{sec:num_spec_eq}, we aim to minimize the discrepancy between $A_{g}v$ and $A_{c}v$ where $v$ is an algebraically smooth vector. 

For each PDE coefficient $\beta \in \Re^p$ possessing a probability distribution $p_{\beta}$ in $\mathscr{B}$, the loss function tied to $F_{\theta}$ and $G_{\psi}$ is defined as:
\begin{equation}
\label{eq:loss}
\mathcal{L}_{\beta}\left(F_{\theta}, G_{\psi}, \mathcal{A}_{g}^\beta, \{v^\beta_{j}\}_{j=1}^{s} 
\right) = 
\sum_{j=1}^{s}
\norm{A_{g}^\beta v^{\beta}_{j} - A_{c}^\beta v^{\beta}_{j}}_{2}^{2},
\end{equation} 
where $\{v^{\beta}_{j}\}$ represents the set of algebraically smooth vectors, 
$s$ is the number of these vectors, 
and $\mathcal{A}_{c}^\beta$ is computed by \cref{alg:sparsify}.
The objective is to minimize the expectation of $\mathcal{L}_{\beta}$ under the distribution of $\beta$, 
symbolized as $\Ee_{\beta \sim p_{\beta}}\left[ \mathcal{L}_{\beta}\right]$, throughout the training.
It is worth mentioning that instead of explicitly forming the matrices $A_{g}^\beta$ and $A_{c}^\beta$, 
we adopt a stencil-based approach where the matrix-vector multiplications are performed 
as the convolutions of the stencils $\mathcal{A}_{g}^\beta$ and $\mathcal{A}_{c}^\beta$ with vectors 
that are padded with zero layers. 
The stencil-based approach and the convolution formulation greatly enhance memory efficiency during training.

\subsection{Details of Training and Testing}
In this section, we provide details of the training and testing 
algorithms.
\paragraph{Architecture of the multi-head attention}  
We use multi-head attention \cite{attention} to compute both location probability in Step 1 with $F_{\theta^{(l)}}$
and numerical values in Step 2 with $G_{\theta^{(l)}}$ 
as discussed in \cref{sec:ML}. We adopt the standard architecture that comprises a set of $n_h$ independent attention heads, each of which extracts different features from each stencil entry of $\mathcal{A}_{g}^{(l)}$. In essence, each head generates a different learned weighted sum of the input values, where the weights are determined by the attention mechanism and reflect the importance of each value. The weights are calculated using a softmax function applied to the scaled dot-product of the input vectors. The output from each head is then concatenated and linearly transformed to produce the final output. 

The multi-head attention mechanism in our study is formally defined as follows:
Let $v_g^{(l)}$ denote the input vectors. For each attention head $h_i,
i=1,\ldots, n_h$, we first transform the inputs using parameterized linear transformations, $W^{Q}_{i}$, $W^{K}_{i}$, and $W^{V}_{i}$ to produce 
the vectors of
query $Q_i$, key $K_i$, and value $V_i$ as follows:
\begin{equation}
Q_{i} = v_g^{(l)} W^{Q}_{i}, \quad
K_{i} = v_g^{(l)} W^{K}_{i}, \quad
V_{i} = v_g^{(l)} W^{V}_{i}.
\end{equation}
The attention scores for each input vector in head $h_i$ are then computed using the scaled dot-product of the query and key vectors, followed by a softmax function:
\begin{equation}
\text{Attention}_{i} = \text{softmax}\left(\frac{Q_{i} K_{i}^T}{\sqrt{d_k}}\right) V_{i} \ ,
\end{equation}
where $d_k$ is the dimension of the key vectors. This process captures the dependencies among the input vectors based on their similarities.
The output of each attention head $h_i$ is then concatenated, and a linear transformation is applied using a parameterized weight matrix $W^{O}$ with softmax activation, which ensures positive outputs:
\begin{equation}
\text{MultiHead}(v_g^{(l)}) = \text{Concat}(\text{Attention}_{1}, \ldots, \text{Attention}_{n_h}) W^{O}.
\end{equation}
This architecture empowers the model to learn complex PDE stencil patterns effectively. The design is flexible, and the number of heads can be adjusted as per the complexity of the task. 

\paragraph{Intuition of selection of multi-head attention}
We first briefly explain why multi-head attention is beneficial for  PDE stencil learning than other types of NNs. This work is about teaching the NNs to generate stencils, which are essentially small patterns or templates used in the discretization of PDEs. These stencils represent the relationship between a grid point and its neighbors.
In the context of PDE stencil learning, multi-head attention can be highly beneficial for several reasons:
\begin{enumerate}[leftmargin=1.5em]
    \item Feature diversification: The multi-head attention allows the model to focus on various features independently, and thus, can capture a wider variety of patterns in the data. For PDE stencil learning, this means that the model can understand the relationships between different grid points more comprehensively.
    \item Context awareness: Attention mechanisms inherently have the capacity to consider the context, i.e., the relationships between different parts of the input data. In PDE stencil learning, this translates to understanding the interactions between a grid point and its surrounding neighboring points.
    \item Flexibility: Multi-head attention adds flexibility to the model. Each head can learn to pay attention to different features, making the model more adaptable. In the context of PDEs, this means that one head can learn to focus on local features (such as the values of nearby points), while another might focus more on global or structural aspects.
\end{enumerate}
We note here that these explanations owe to the empirical observation that such an architecture works better than  vanilla deep NNs on our task.
\paragraph{{Details of training and testing}}
The PDE coefficient $\beta_i$ is sampled 
from distribution $\mathscr{B}$ according to the probability 
density function
$p_{\beta}$ to get the set of $N_t$ parameters, $\{\beta_i\}$, $i=1, \ldots, {N_t}$. 
Then, we construct the corresponding set of fine-grid stencils $\{\mathcal{A}^{\beta_i}\}$.
%\textcolor{red}{H: please improve, An important aspect to consider is that while we utilize multi-head attention 
%to learn the locations and values of nonzero entries, we also ensure that the output of our neural network has mathematical significance. To accomplish this, we intentionally avoid learning the 
%central point of $\mathcal{A}_c$. 
%Subsequently, we set the value of the central point of $\mathcal{A}_c$ to be the negative sum of the learned outputs. }
For all the tests in this paper, we use full coarsening, full-weighting restriction,
and the corresponding bi-linear interpolation for all the levels of AMG.
At each level $l>1$, the ML model is built with
the Galerkin stencils
$\{(\mathcal{A}_{g}^{\beta_i})^{(l)}\}$
and a set of smooth test vectors 
$\{(v_{j}^{\beta_i})^{(l)}\}$, $j=1, \ldots, {s_l}$, associated with each of the stencil,
using
the loss function 
\begin{equation}\label{eq:lossL}
\mathcal{L}
%\left(F_{\theta^{(l)}}, G_{\psi^{(l)}}, \{(\mathcal{A}_g^{\beta_i})^{(l)}\}%_{i=1}^{N_t}
%, \{(v_{j}^{\beta_i})^{(l)}\}\right) 
=
\frac{1}{N_t}\sum_{i=1}^{N_t}\mathcal{L}_{\beta_i}\left(F_{\theta^{(l)}}, G_{\psi^{(l)}}, (\mathcal{A}_g^{\beta_i})^{(l)}, \{(v^{\beta_i}_{j})^{(l)}\} \right)
\end{equation}
that
is used to approximate
$\Ee_{\beta \sim p_{\beta}}\left[ \mathcal{L}_{\beta}\right]$.
The complete training procedure is summarized in \cref{alg:train}. 
It is important to note that the NN trainings are independent of each other on different levels. 
Therefore, the training of the NNs for each level can be carried 
out simultaneously once the training data is prepared, taking advantage of parallel computing resources.
The testing set is constructed with parameters that differ from those in the 
training set. This means that we test the model on a set of PDE parameters $\{\beta_j\}$, $j=1, \ldots, {N_v}$, 
that have not been encountered by the models during training. 
%Here, $N_v$ represents the size of our testing set. 
The purpose of the testing set 
is to assess the generalization capability of new problem instances.

\begin{algorithm}[h]
\caption{Training NNs for computing coarse-grid operator at level $l$}
\label{alg:train}
\begin{algorithmic}[1]

\REQUIRE Interpolation operator $P^{(l)}$, 
Galerkin coarse-grid stencils $\{(\mathcal{A}_g^{\beta_i})^{(l)}\}$,  
the number of test vectors $s_l$ and the target stencil complexity $k$

\ENSURE NNs $F_{\theta^{(l)}}$ and $G_{\psi^{(l)}}$ 

 \STATE Generate test vectors $(v^{\beta_i}_{j})^{(l)}$, $j=1,\ldots,s_l$, for each $(\mathcal{A}_g^{\beta_i})^{(l)}$ as follows %using \cref{alg:smooth-basis}
\begin{itemize} 
\setlength{\itemindent}{-2em}
\item Compute $T^{(l)}=(P^{(l)})\trans P^{(l)}$
\item Compute the eigenvalues and vectors of 
\begin{equation} \notag
(A_g^{\beta_i})^{(l)} u = \lambda_i^{(l)} T^{(l)} u
\end{equation}
\item The test vectors are the eigenvectors associated with the $s_l$ smallest eigenvalues
\end{itemize}
 \STATE Initialize $G_{\theta^{(l)}}$ and $G_{\psi^{(l)}}$
\REPEAT
\STATE 
Apply \cref{alg:sparsify}:
$(\mathcal{A}_{c}^{\beta_i})^{(l)}=
\mbox{SparsifyStencil}
\left((\mathcal{A}_{g}^{\beta_i})^{(l)}, F_{\theta^{(l)}}, G_{\psi^{(l)}}, k\right)$
 \STATE Compute the gradient of the loss \cref{eq:lossL}
%\begin{equation*}
%$
%    \mathcal{L} = \frac{1}{N_t}\textstyle\sum_{i=1}^{N_t} 
%    \sum_{j=1}^{s_l}||A_{\beta_i}^{(l)}v^{(l)}_{\beta_i,j} -
%A_{\beta_i,nn}^{(l)}v^{(l)}_{\beta_i,j}||_{2}^{2}
%$
%\end{equation*}
\STATE Update the weights $\theta^{(l)}$ and $\psi^{(l)}$
\UNTIL{the prescribed number of training epochs is reached}
\end{algorithmic}
\end{algorithm}

% \begin{algorithm}
% \caption{Algebraically smooth test vectors}
% \label{alg:smooth-basis}
% \begin{algorithmic}[1]
% %\REQUIRE Interpolation operator $P$, coarse-level stencil $\mathcal{A}_{g}$, number 
% %of algebraically smooth vectors: $s$, size of the eigenvalue problem: $n$
% %\ENSURE Algebraically smooth basis $\bfv_{i}, i=1,\dots,k$ 
% \STATE Compute $T=P\trans P$
% \STATE Compute the generalized eigenvalues and vectors of $A_{g}v_i=\lambda_i T v_i$.
% \RETURN $\{v_{i}\}$ associated with the $k$ smallest eigenvalues
% \end{algorithmic}
% \end{algorithm}

\section{Numerical Results}
\label{sec-numerical-exp}

We report the numerical results of the proposed ML-based
non-Galerkin coarse-grid method in this section. 
All the ML models in the work\footnote{The codes is available at \url{https://anonymous.4open.science/r/Sparse-Coarse-Operator-11C7}} were written with
PyTorch 1.9.0 \cite{paszke2019pytorch}. 
We use PyAMG 4.2.3 \cite{bell2022pyamg} to build the AMG hierarchy. 
All the experiments were performed on a workstation with Intel Core i7-6700 CPUs.

\subsection{Evaluation Metrics}
In this section, we evaluate the performance of the proposed ML-based approach  
by comparing the average number of iterations required 
by the AMG method  
using different coarse-grid operators
to converge.
%to a predefined tolerance of the residual norm. 
Additionally, we verify 
the spectral equivalence of the Galerkin and sparsified non-Galerkin stencils
by computing the spectra of the corresponding matrices on  
meshes of various sizes.
\subsection{Spectrally equivalent stencils}
%\textcolor{red}{RL: I was here, as of 6/24 13:25.}
We first examine the proposed ML-based method on the 9-point stencil problem
\cref{eq:9ptcirculant}
that allows direct evaluation of the learned non-Galerkin operator 
by the comparison with 
the theoretical result
\cref{eq:5ptcirculant}, which is a spectrally equivalent
5-point stencil.
% \textbf{***TODO:***}
% \begin{todolist}
%     \item introduce generating symbol? 
% \end{todolist}
We use the 9-point stencil $\mathcal{A}$ of the form \cref{eq:9ptcirculant} with
$a = 2.720$, $b = 1.417$ and $c = 0.000114$, i.e.,
\begin{align}
   \mathcal{A} 
       &=
       \left[
       \begin{array}{ccc}
       0.000114&2.720&0.000114\\
       1.417&-8.275&1.417\\
       0.000114&2.720&0.000114\\
       \end{array}\label{eq:circ-stencil}
       \right]
\end{align}
as the fine-level $A$-operator, and
the 2-D full-weighting stencil,
\begin{equation}
\mathcal{R} = \dfrac{1}{4}\begin{bmatrix}
1 & 2 & 1 \\
2 & 4 & 2 \\
1 & 2 & 1 
\end{bmatrix},
\end{equation}
for the restriction operator.
Thus, the stencil of the Galerkin operator is
\begin{equation} \label{eq:Ag}
\mathcal{A}_g = \left[
        \begin{array}{ccc}
        0.129 & 0.096 & 0.129\\
        0.422 & -1.551 & 0.422\\
        0.129 & 0.096 & 0.129\\
        \end{array}
        \right],
\end{equation}
which has the same form as $\mathcal{A}$.
From \cref{eq:5ptcirculant}, a 5-point stencil that is 
spectrally equivalent to \cref{eq:Ag} is given by 
\begin{equation} \label{eq:Ac}
\mathcal{A}_{c} = \begin{bmatrix}
    0 & 0.354 & 0 \\
    0.680 & -2.069 & 0.680 \\
    0 & 0.354 & 0 
    \end{bmatrix}.
\end{equation}
Using \cref{alg:train} with the prescribed stencil complexity
$k=5$ on the $31\times 31$ grid, 
the pre-trained NNs produced the following 5-point stencil
\begin{equation}
\mathcal{A}_{nn} = \begin{bmatrix}
    0 & 0.348 & 0 \\
    0.666 & -2.024 & 0.663 \\
    0 & 0.347 & 0 
    \end{bmatrix},
\end{equation}
denoted by $\mathcal{A}_{nn}$, which is 
close to the theoretical result \cref{eq:Ac}.
To assess the convergence behavior of the AMG method, 
we solve a linear system using the coefficient matrix defined by \eqref{eq:circ-stencil}.
We conduct these tests on larger-sized grids and use
the two-grid AMG methods employing 
$\mathcal{A}_g$, $ \mathcal{A}_{c}$, and $\mathcal{A}_{nn}$
as the coarse-level operators, respectively.
The right-hand-side vector is generated randomly.
The stopping criterion with respect to the relative residual
norm is set to be $10^{-6}$. The results are shown
in  \cref{tbl:circ-result}, 
from which we can observe that all three methods require the same number of iterations.

\begin{table}[h!t]
    \centering
    \caption{The number of iterations required by 
    the two-level AMG methods for solving 
    a linear system corresponding to the coefficient matrix stencil \eqref{eq:circ-stencil}
    to $10^{-6}$ accuracy in terms
    of the relative residual.
    The AMG methods utilize $\mathcal{A}_g$, $ \mathcal{A}_c$, and $\mathcal{A}_{nn}$ 
    as the coarse-level operator respectively.
    The tests are carried out on grid sizes up to $511 \times 511$.}
    \begin{tabular}{c|c|c|c|c|c|c|c}
    & \multicolumn{7}{c}{grid size } \tabularnewline
    \cline{2-8}
    & 63 & 95 & 127 & 191 & 255 & 383 & 511\\\hline
    {$\mathcal{A}_{g}$} & 11  & 10  & 10 & 10  & 10 & 10  & 10 \\
    {$\mathcal{A}_{c}$}& 11 & 10 & 10 & 10  & 10 & 10  & 10\\
    {$\mathcal{A}_{nn}$} & 11 & 10 & 10 & 10  & 10 & 10  & 10\\
    \end{tabular}
    \label{tbl:circ-result}
\end{table}

\subsection{2-D rotated Laplacian problem}
In this section, we consider the 2-D anisotropic rotated Laplacian problem
\begin{equation}\label{pde}
    -\nabla \cdot (T_{\theta,\xi}\nabla u(x,y)) = f(x,y), %\quad  (x,y) \in \mathcal{G}
\end{equation}
where the $2\times 2$ vector field $T_{\theta,\xi}$ parameterized by $\theta$ and $\xi$ is defined as
\begin{equation}
\notag
T_{\theta,\xi}=\begin{bmatrix}
\cos^{2}{\theta}+\xi\sin^{2}{\theta}&\cos{\theta}\sin{\theta}(1-\xi) \\ 
\cos{\theta}\sin{\theta}(1-\xi) & \sin^{2}{\theta}+\xi\cos^{2}{\theta}
\end{bmatrix}
\end{equation}
with $\theta$ being the angle of the anisotropy and $\xi$ being the conductivity.

We show that the proposed approach is not limited to a particular set of parameters 
but remains effective across a range of values for both $\xi$ and $\theta$.
In the first set of experiments, we fix the value $\xi$ while allowing $\theta$ to follow a uniform distribution within a specified interval. We conduct 12  experiments 
where each $\xi \in \{100,200,300,400\}$ is paired up with $\theta$ sampled from 
intervals $\{({\pi}/{4}, {\pi}/{3}), ({\pi}/{3}, {5\pi}/{12}), 
({\pi}/{2},{7\pi}/{12})\}$. 
The AMG methods for solving these problems use
full-coarsening, full-weighting restriction, and 
the Gauss-Seidel method for both pre-smoothing and post-smoothing.
AMG V-cycles are executed until the residual norm is reduced by $6$ orders of magnitude. 
The number of the nonzero elements is $9$ in the Galerkin coarse-grid stencil across
the AMG levels, whereas
we choose to reduce the number to $5$ for the non-Galerkin operator. 

During the training phase of each experiment, the model is provided with 5 distinct 
instances that share the same  $\xi$ but have different $\theta$ values evenly distributed within the chosen interval. 
For example, for $\xi=100$ and $\theta \in ({\pi}/{4}, {\pi}/{3})$, the parameters for the $5$ instances of $(\xi,\theta)$ are selected as follows:
\begin{equation}
\left\{(100, {\pi}/{4}),(100, {3.25\pi}/{12}), 
(100, {3.5\pi}/{12}),(100, {3.75\pi}/{12}),(100, {\pi}/{3})\right\}.
\end{equation}
The size of the fine-level matrix in the training instances is set to be $31\times 31$.
In the testing phase, $10$ distinct $\theta$ values are randomly selected from the chosen
interval. The AMG parameters are consistent with those used in the training phase.
In the testing, it should be noted that the fine-level problem size is 
$511\times 511$, which is approximately 256 times larger than that in the training instances. This larger problem size in the testing allows for a more rigorous evaluation of the performance of AMG and the ability to handle larger-scale problems.
We record the number of iterations required by the 3-level AMG method to converge with
the Galerkin and non-Galerkin operators, 
shown in  \cref{tbl:laplacian-vary-theta}.
These results indicate that the convergence behavior of the AMG method remains largely unchanged when the alternative sparser non-Galerkin coarse-grid operators are used as replacements.

\begin{table}[h!]
\centering
\caption[Averaged numbers of iterations]{The average number of iterations 
required by the 3-level AMG to converge with the Galerkin and non-Galerkin 
coarse-grid operators 
for solving \eqref{pde} with different PDE parameter $\xi$ and $\theta$. 
The mesh size is $511\times 511$. 
The parameters are selected so that $\xi \in \{100,200,300,400\}$ is fixed and $\theta$ is 
randomly sampled from a uniform distribution in each interval. 
The iteration number is averaged over 10 different sampled $\theta$ values.}

\begin{tabular}{c|c|c|c|c  }
 %\hline
 & \multirow{2}{*}{$\xi$} & \multicolumn{3}{c}{$\theta$}  \\
\cline{3-5}
 &  & (${\pi}/{6}$, ${\pi}/{4}$) & (${\pi}/{4}$,${\pi}/{3}$) & 
 (${\pi}/{2}$, ${7\pi}/{12}$) \\
 \hline
 $\mathcal{A}_{g}$  & \multirow{2}{*}{100} &92.1 & 102.8 &126.9\\
 ${\mathcal{A}_{nn}}$ &  &89.0 & 93.0 &135.2\\
 \hline
 \hline
 $\mathcal{A}_{g}$ & \multirow{2}{*}{200} & 191.7 &196.6 & 203.1\\
 ${\mathcal{A}_{nn}}$ & &174.2 &177.8 & 204.9\\
 \hline
 \hline
 $\mathcal{A}_{g}$ & \multirow{2}{*}{300} & 248.0 &269.7 & 342.3\\
 ${\mathcal{A}_{nn}}$ & & 246.5 &231.4 & 356.2\\
 \hline
  \hline
 $\mathcal{A}_{g}$ & \multirow{2}{*}{400} & 337.1 &351.1 & 438.2 \\
 ${\mathcal{A}_{nn}}$ & & 326.3 &327.7 & 441.5 \\
 \hline
\end{tabular}
\label{tbl:laplacian-vary-theta}
\end{table}

% For evaluating the metric $\phi$ in \eqref{eq:phi}, we take a specific problem with $\xi=10$ and $\theta=\frac{\pi}{3}$ as an example. We show the eigenvalues of $\mxA_{c}^{-1}\mxA_{g}$ for several different mesh sizes in Figure \ref{fig:eigenvalues-ac-ag}. Note that for apparent practical issues, we may only present for several different mesh sizes. We can see that the distribution of spectrum is the same for the selected sizes, empirically indicating spectral equivalence of $\mxA_{c}$ and $\mxA_{g}$.

% \begin{figure}[]
% \centering
% \includegraphics[width=0.9\columnwidth]{figs/eig_value_diffusion_1.png}
% \caption[Eigenvalues]{The eigenvalues of $\mxA_c^{-1}\mxA_{g}$ with different matrix sizes for solving rotated Laplacian problem with $\xi=10$ and $\theta=\frac{\pi}{3}$.}
% \label{fig:eigenvalues-ac-ag}
% \end{figure}

In the second set of experiments, we keep the parameter $\theta$ fixed
and vary $\xi$ following a uniform
distribution within the selected intervals. 
A total of 12 experiments were conducted
where each $\theta \in \{{\pi}/{6},{\pi}/{4},{\pi}/{3},{5\pi}/{12}\}$ is paired with $\xi$ sampled from the intervals $\{(5,10), (80,100), (100,200)\}$. 
The AMG configurations used in these experiments remain the same as in the previous set.
The training and testing processes are also similar. 
For each experiment, we train the model 
using 5 different instances evenly distributed 
within the selected intervals and then test it with 10 
randomly sampled $\xi$ values from the same interval.
The size of the fine-level linear system in the training instances is set to be $31\times 31$, while in each testing instance, it has a much larger size that is $511\times 511$. The averaged numbers of iterations from all the experiments are presented in  \cref{tbl:laplacian-vary-xi}.

\begin{table}[htb]
\centering
\caption[Averaged numbers of iterations]{The average number of iterations 
required by the 3-level AMG to converge with the Galerkin and non-Galerkin 
coarse-grid operators 
for solving \eqref{pde} with different PDE parameter $\xi$ and $\theta$. 
The mesh size is $511\times 511$. 
The parameters are selected such that $\theta \in 
\{{\pi}/{6},{\pi}/{4},{\pi}/{3},{5\pi}/{12}\}$ is fixed and $\xi$ is randomly sampled 
from a uniform distribution in each interval. 
The iteration number is averaged over  10 different sampled $\xi$ values.}

% \begin{tabular}{c|c|c|c}
% $\theta=\frac{\pi}{4},\xi$ &(0.1,0.13) &(0.13,0.16) &(0.16,0.19)  \\\hline
% $\mathcal{A}_{g}$ &13.1  &11.4 &10.2 \\
% ${\mathcal{A}_{c}}$&17.7 &11.8&10.4 \\

% \end{tabular}

\begin{tabular}{c|c|c|c|c  }
 %\hline
 & \multirow{2}{*}{$\theta$} & \multicolumn{3}{c}{$\xi$}  \\
\cline{3-5}
 &  & (100, 200) & (80, 100) & 
 (5, 10) \\
 \hline
 $\mathcal{A}_{g}$  & \multirow{2}{*}{$\pi/6$} & 90.4 & 72.1 &13.5\\
 ${\mathcal{A}_{nn}}$ &  &100.2 & 84.4 &13.8\\
 \hline
 \hline
 $\mathcal{A}_{g}$ & \multirow{2}{*}{$\pi/4$} & 172.5 &105.2 & 14.1\\
 ${\mathcal{A}_{nn}}$ & &123.1 &79.0 & 15.9\\
 \hline
 \hline
 $\mathcal{A}_{g}$ & \multirow{2}{*}{$\pi/3$} & 99.4 & 80.9 & 14.3\\
 ${\mathcal{A}_{nn}}$ & & 79.1 & 88.8 & 15.4\\
 \hline
  \hline
 $\mathcal{A}_{g}$ & \multirow{2}{*}{$5\pi/12$} & 92.5 & 76.4 & 16.5 \\
 ${\mathcal{A}_{nn}}$ & & 107.4 & 88.2 & 16.6 \\
 \hline
\end{tabular}
\label{tbl:laplacian-vary-xi}
\end{table}

% \noindent We also present the comparison between the solutions using dense and learned sparsified stencils for a sampled test case in Figure \ref{fig:solution-ac-ag} to give a straightforward evaluation on the performance of our algorithm. 

% \begin{figure}[H]
% \centering
% \label{fig:solution-ac-ag}
% \includegraphics[width=0.45\columnwidth]{figs/rotated_solution_ag.png}
% \includegraphics[width=0.45\columnwidth]{figs/rotated_solution_ac.png}
% \caption[The approximate solution]{The approximate solutions (left: $\mathcal{A}_{g}$, right: ${\mathcal{A}_{c}}$) after 10 iterations of two-grid method with two coarse grid operators for solving rotated Laplacian problem with $\xi=0.1147$ and $\theta=\frac{\pi}{4}$.}
% \end{figure}

\begin{figure}[h]
\centering
\includegraphics[width=\columnwidth]{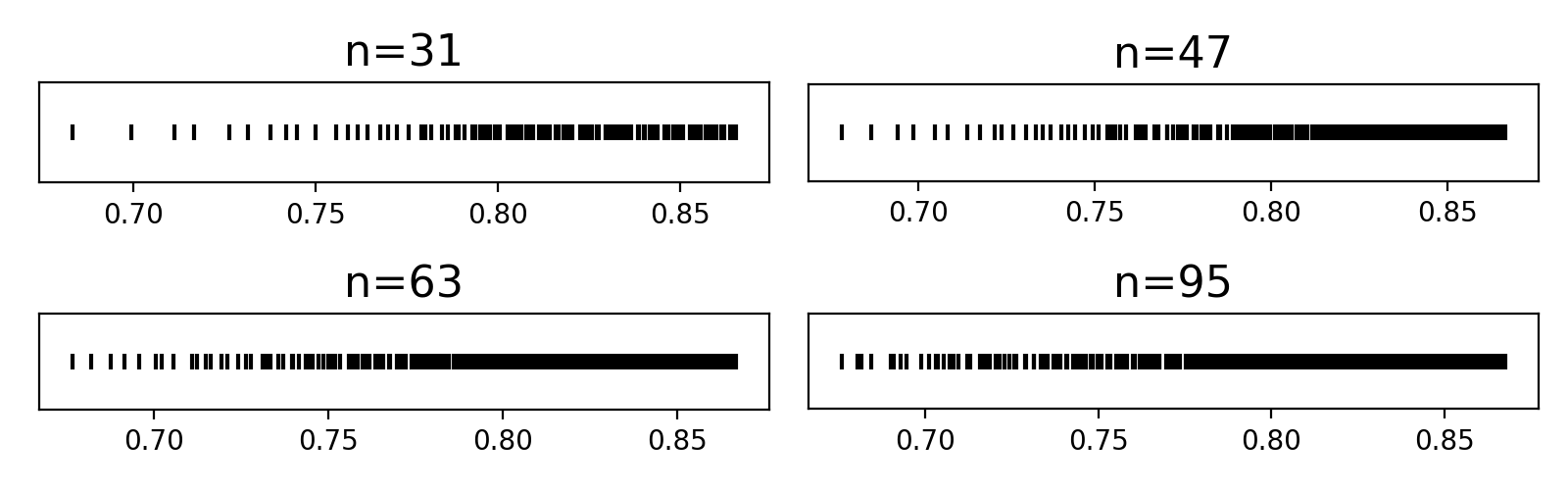}
\caption[Eigenvalues]{The eigenvalues of $A_{nn}^{-1}A_{g}$ 
on meshes of different sizes ($n\times n$) for solving the rotated Laplacian problem with $\theta={\pi}/{6}$ and $\xi=0.1$.}
\label{fig:eigenvalues-ac-ag-elas}
\end{figure}

In the subsequent experiment, we specifically consider the Laplacian problem with parameters
$\theta={\pi}/{6}$ and $\xi=0.1$ 
as an example to demonstrate the measurement of spectral equivalence
as defined in \cref{eq:speq}.
We examine the eigenvalues of $A_{nn}^{-1}A_{g}$ on meshes of varying sizes, as depicted in \cref{fig:eigenvalues-ac-ag-elas}. 
We observe that all the eigenvalues are approximately equal to 1, and the distribution of eigenvalues remains consistent regardless of the mesh size. 
This observation suggests the presence of spectral equivalence between the two coarse-grid operators across meshes of different sizes. 

\begin{figure}[h]
\centerline{\includegraphics[width=0.5\columnwidth]{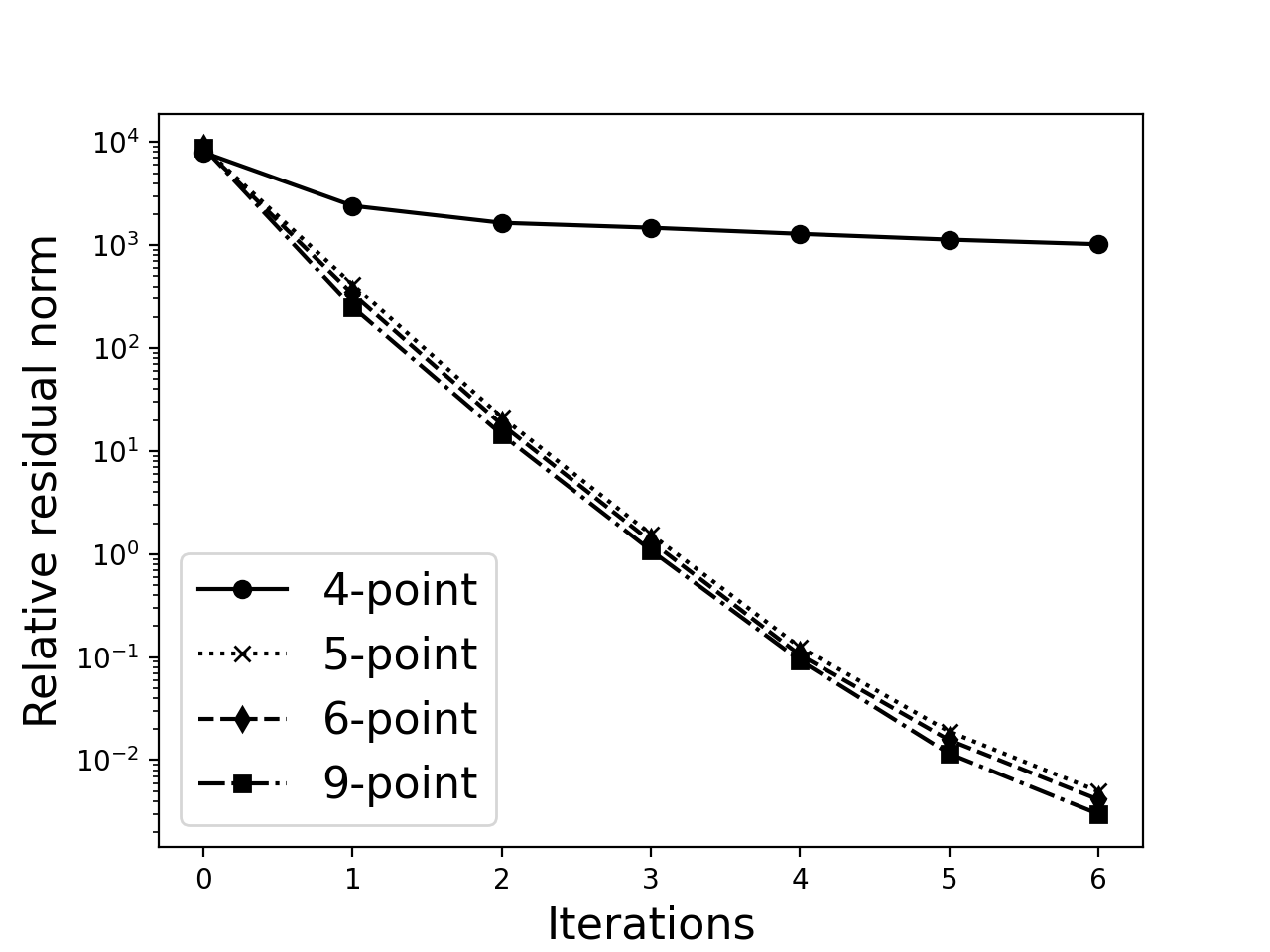}}
\caption[Convergence]{The convergence in terms of
the residual norm 
of the two-grid AMG methods using the coarse-grid operator
from the NN model of stencil complexity $k=4, 5,6$ and the 
Galerkin operator
for solving the rotated Laplacian problem with $\xi=10$ and $\theta={\pi}/{4}$.}
\label{fig:convergence-ac-ag}
\end{figure}

The target stencil complexity $k$
in \cref{alg:sparsify} is a parameter 
left to be chosen by the users.
It is an adjustable parameter that allows users to control the sparsity level in
the trained NN-model
and of the resulting coarse-grid operator. The appropriate value of $k$ typically depends on the problem domain and the desired balance between accuracy and computational efficiency.
It may be necessary to perform experiments 
to determine the optimal value of $k$ for a particular 
application.
In the final experiment, we perform this study
for the rotated Laplacian problem with $\xi=10$ and 
$\theta = {\pi}/{4}$.
Note that the Galerkin operator has a 9-point stencil,
so we vary the stencil complexity from $4$ to $6$ 
in the non-Galerkin operator
and 
record the convergence behavior of the corresponding AMG method.
The results, as depicted in  \cref{fig:convergence-ac-ag}, show the findings regarding the convergence behavior of the AMG method with different values of $k$ in the sparsified stencil  $\mathcal{A}_{nn}$.
Notably, when $k=4$, the AMG method fails to converge. However, for $k=5$ and $k=6$, the convergence behavior closely resembles that of the 9-point Galerkin operator. This observation suggests that a minimum stencil complexity of $k=5$ appears to be required for $\mathcal{A}_{nn}$ to achieve convergence, which coincides with the operator complexity of the fine-grid operator.

\subsection{2-D linear elasticity problem}\label{sec:elasticity}

In this section, we consider the 2-D time-independent linear elasticity problem in an isotropic homogeneous medium:
\begin{align} \label{eq:elast}
    \mu \nabla^{2}u+(\mu+\lambda)(\frac{\partial^{2}u}{\partial x^{2}}+\frac{\partial^{2}v}{\partial x\partial y})+f_{x} &= 0,\\
    \mu \nabla^{2}v+(\mu+\lambda)(\frac{\partial^{2}v}{\partial x^{2}}+\frac{\partial^{2}u}{\partial x\partial y})+f_{y} &= 0,
\end{align}
where $u$ and $v$ are the solution in the direction of $x$- and $y$-axis respectively, $f_{x}$ and $f_{y}$ are forcing terms, and
$\mu$ and $\lambda$ are Lame coefficients that are determined by Young's modulus $E$ and Poisson's ratio $\nu$ as
\begin{equation}
\mu = \frac{E}{2(1+\nu)}\ , \quad \lambda = \frac{E\nu}{(1+\nu)(1-2\nu)} \ .
\end{equation}
In our tests, we set $E=10^{-5}$ and vary the value of $\nu$.
For the discretization, 
we adopt the optimal 2-D 9-point stencil
in terms of local truncation errors \cite{idesman2020compact} on rectangular
Cartesian grid with the mesh step sizes
$h$ and $b_{y}h$, respectively,  
along the $x$- and $y$-axes ($b_y$ is the aspect ratio of the mesh),
\begin{align}
\mathcal{A}_{uu} &=
\begin{bmatrix}
a_{uu}^{nw}  &  a_{uu}^{n}  & a_{uu}^{ne}\\
a_{uu}^{w} &1  & a_{uu}^{e}  \\
a_{uu}^{sw} & a_{uu}^{s} & a_{uu}^{se} \\
\end{bmatrix}, \quad
\mathcal{A}_{uv} =
\begin{bmatrix}
 a_{uv}^{nw} &0  &  a_{uv}^{ne} \\
0 & 0 &  0 \\
a_{uv}^{sw} & 0 & a_{uv}^{se} \\
\end{bmatrix}, \;
\end{align}
where the coefficients are given by
\begin{align}
a_{uu}^{n} &= a_{uu}^{s} = \frac{(b_{y}^{2}-1)\lambda+2(b_{y}^{2}-2)\mu}{2(2\lambda b_{y}^{2}+\lambda+4(b_{y}^{2}+1)\mu)} \ , \\
a_{uu}^{w} &= a_{uu}^{e} = -\frac{2\lambda b_{y}^{2}+4\mu b_{y}^{2}+\lambda+\mu}{2(2\lambda b_{y}^{2}+\lambda+4(b_{y}^{2}+1)\mu)} \ , \\
a_{uu}^{nw} &= a_{uu}^{ne} = a_{uu}^{sw} = a_{uu}^{se} = \frac{-\lambda b_{y}^{2}-2\mu b_{y}^{2}+\lambda+\mu}{4(2\lambda b_{y}^{2}+\lambda+4(b_{y}^{2}+1)\mu)} \ , \\
a_{uv}^{nw} &= a_{uv}^{se} = -a_{uv}^{ne} = -a_{uv}^{sw} = \frac{3b_{y}(\lambda+\mu)}{8(2\lambda b_{y}^{2}+\lambda+4(b_{y}^{2}+1)\mu)} \ .
\end{align}
These stencils define the $2 \times 2$ block linear system 
%corresponding to \cref{eq:elast}
\begin{equation}\label{eq:linear-system-elasticity}
\begin{bmatrix}A_{uu}& A_{uv}\\ A_{vu} &A_{vv} \end{bmatrix}
\begin{bmatrix}u\\ v\end{bmatrix}=\begin{bmatrix} f_{x}\\ f_{y}\end{bmatrix},
\end{equation}
where ${A}_{uu} = {A}_{vv}$ and ${A}_{vu} = {A}_{uv}\trans$.
A node-based AMG approach is used to solve
\cref{eq:linear-system-elasticity}, where the same
red-black coarsening is used in
$u\mhyphen u$ and $v \mhyphen v$ blocks and the interpolation and restriction
operators have the same block form
\begin{equation} \label{eq:interp_restric}
R=
\begin{bmatrix}
R_{uu} & R_{uv}\\
R_{vu} & R_{vv}
\end{bmatrix},
\quad 
P=
\begin{bmatrix} 
P_{uu} & P_{uv}\\
P_{vu} & P_{vv} 
\end{bmatrix} ,
\end{equation}
which interpolate and restrict within and across the 
two types of variables $u$ and $v$.
The stencils of the 
operators in \cref{eq:interp_restric} are given by, respectively,
\begin{align}
\mathcal{R}_{uu} &= \mathcal{R}_{vv} =
\frac{1}{8}
\left[
\arraycolsep=2pt
\begin{array}{ccc}
&1&\\
1&4&1\\
&1&\\
\end{array}
\right] ,
\quad
\mathcal{P}_{uu} = \mathcal{P}_{vv} =
\frac{1}{4}
\left]
\arraycolsep=2pt
\begin{array}{ccc}
&1&\\
1&4&1\\
&1&\\
\end{array}
\right[ \ ,
\\
\mathcal{R}_{uv} &=
\frac{1}{8}
\left[
\arraycolsep=2pt
\begin{array}{ccc}
&1&\\
-1&0&-1\\
&1&\\
\end{array}
\right] ,
\quad
\mathcal{P}_{uv} = 
\frac{1}{4}
\left]
\arraycolsep=2pt
\begin{array}{ccc}
&1&\\
-1&0&-1\\
&1&\\
\end{array}
\right[ \ ,
\\
\mathcal{R}_{vu} &=
\frac{1}{8}
\left[
\arraycolsep=2pt
\begin{array}{ccc}
&-1&\\
1&0&1\\
&-1&\\
\end{array}
\right],
\quad 
\mathcal{P}_{vu} =
\frac{1}{4}
\left]
\arraycolsep=2pt
\begin{array}{ccc}
&-1&\\
1&0&1\\
&-1&\\
\end{array}
\right[ \ .
\end{align}
As stated in \cite{brezina2001algebraic},
to interpolate exactly the smoothest function that is locally constant, 
it requires the
interpolation weights
for $u\mhyphen u$ and $v\mhyphen v$ to sum to 1 and 
for the $u\mhyphen v$ and $v\mhyphen u$ weights to sum to 0.
The Gauss-Seidel smoother is used with the AMG V-cycle
and the iterations are stopped when the
relative residual norm is below $10^{-6}$. 

We train the NN model on 4 different instances with 
$\nu \in \{0.1,0.2,0.3,0.4\}$ to
reduce the complexity of the Galerkin operator by 50\%.
The coarse-grid Galerkin operator has the same block structure
as \cref{eq:linear-system-elasticity} and only
2 distinct stencils due to the symmetry of the matrix.
In the training, we combine these 2 stencils and pass them to
the NNs as the input. 
It turns out that the NN model trained in this way yields
better coarse-grid operators than learning the stencils of
$u\mhyphen u$ and $u\mhyphen v$ separately.  

The mesh size used in the training set is set to be $9\times 9$.
We test the model on instances with 
$\nu$ randomly drawn from 
each interval of $\{(0.1,0.2), (0.2,0.3), (0.3,0.4)\}$. 
The size of the mesh used in the testing is 
$65\times 65$.
The average numbers of iterations are presented in \cref{tbl:sp-lin-elas}. Similar to the results observed in the rotated Laplacian problems, 
the convergence behavior of the two-grid AMG method is not negatively affected by the replacement with the non-Galerkin coarse-grid operator obtained from the NN model.

\begin{table}[h!]
\centering
\caption{The average number of iterations 
required by the 2-grid AMG to converge with the Galerkin and non-Galerkin 
coarse-grid operators 
for solving \eqref{eq:elast} with 10 different Poisson's ratios
$\nu$ randomly sampled from each interval. The mesh size is set to be $65\times 65$.}

\begin{tabular}{c|c|c|c}
$\nu$& $(0.1,0.2)$ & $(0.2,0.3)$ & $(0.3,0.4)$ \\\hline
$\mathcal{A}_{g}$ & 10.1  & 10.2  & 10.6 \\
${\mathcal{A}_{nn}}$& 11.0 & 10.7 & 11.5 \\
\end{tabular}
\label{tbl:sp-lin-elas}
\end{table}

% The values $\phi = \norm{I-A_{c}A_{g}^{-1}}_{2}$ (\ref{eq:spec-equiv}) for each problem are shown in Table \ref{tbl:sp-lin-elas-phi}. All of the values are less than 1. 

%     \begin{table}[H]
%         \centering
%         \begin{tabular}{c|c|c|c}
%         $\nu$& (0.1,0.2) & (0.2,0.3) & (0.3,0.4) \\\hline
%        $\phi$ &0.7014  & 0.6889 & 0.6890 \\
%         \end{tabular}
%         \caption[$\phi$ for elasticity]{The value $\phi$ (\ref{eq:spec-equiv}) for
%  $E=10^{-5}$ and random $\nu$ on each interval.}
%  \label{tbl:sp-lin-elas-phi}
%     \end{table}

% \noindent An example of the approximate solutions after the same number of iterations of using two-grid method with $\mathcal{A}_{g}$ and ${\mathcal{A}_{c}}$ as coarse grid stencils is presented in Figure \ref{fig:solution-ac-ag-lin-elas} for direct evaluation. 

% \begin{figure}[H]
% \centering
% \includegraphics[width=0.45\columnwidth]{figs/u_solution_ag.png}
% \includegraphics[width=0.45\columnwidth]{figs/v_solution_ag.vpng.png}

% \includegraphics[width=0.45\columnwidth]{figs/u_solution_ac.png}
% \includegraphics[width=0.45\columnwidth]{figs/v_solution_ac.png}

% \caption[The approximate solution]{The approximate solutions to \eqref{eq:linear-system-elasticity} with $E=10^{-5}$ and $\nu=0.1409$ (top left: $\bfu$ solved with $\mathcal{A}_g$, top right: $\bfv$ solved with $\mathcal{A}_g$, bottom left: $\bfu$ solved with $\mathcal{A}_c$, bottom right: $\bfv$ solved with $\mathcal{A}_c$) after 10 iterations of the two-grid scheme.}
% \label{fig:solution-ac-ag-lin-elas}
% \end{figure}
\subsection{Comparison with existing non-Galerkin methods}
In this section, we compare the performance of 
the proposed NN-based algorithm with the  Sparsified Smooth Aggregation
(SpSA) method proposed in \cite{treister2015non} for solving 
the rotated Laplacian problem. 
The SpSA method is based on Smooth Aggregation (SA) AMG methods.
In these methods, a tentative aggregation-based interpolation operator
$P_t$ is first constructed, followed by a few steps of smoothing
of $P_t$ that generate the final interpolation operator $P$,
which is typically considerably denser than $P_t$.
The SpSA algorithm aims to reduce the complexity of the Galerkin 
operator $P\trans A P$ to have the same sparsity pattern as 
$P_t\trans A P_t$.
Given that we utilize the standard Ruge-St\"uben AMG (as opposed to SA AMG) combined with the NN-based approach, conducting a direct comparison 
between the two approaches
becomes challenging due to
the different
AMG hierarchies obtained.
To ensure an equitable comparison, we impose a requirement that the number of nonzero entries per row in the coarse-level operator generated by SpSA should not be smaller than the operator produced by our algorithm.
Consequently, any observed disparities in performance can be attributed to the specific characteristics of the 
selected sparsity pattern and numerical values
of the coarse-grid operator, rather than the variations in the level of the sparsity.
The number of iterations required by 
the GMRES method preconditioned by
3-level AMG methods
with different coarse-grid operators 
for solving the rotated Laplacian problem \eqref{pde} 
are presented in \cref{tbl:laplacian-vary-theta-spsa}
and \cref{tbl:laplacian-vary-xi-spsa}, with
varied PDE coefficients.
For the majority of cases, the AMG method with NN-based coarse-grid operators exhibits better performance compared to SpSA, as it requires fewer iterations to converge to the $10^{-6}$ stopping tolerance and achieves convergence rate that is
much closer to that using the Galerkin operator.
There are a few exceptions where SpSA outperforms the NN-based method, and in some cases, it performs even better than the AMG method using the Galerkin operator.
%We set the relative tolerance for both methods as $10^{-6}$.

\begin{table}[h]
\centering
\caption{The average number of iterations required by 
the GMRES method preconditioned by
3-level AMG methods
with different coarse-grid operators 
for solving \eqref{pde} with 
different sets of PDE parameters. 
The mesh size is $511\times 511$. 
The parameters are selected so that 
$\theta \in \{{\pi}/{6}, {\pi}/{4}, {\pi}/{3},
{5\pi}/{12}\}$
is fixed and $\xi$ is randomly sampled from a uniform distribution in each interval. 
The iteration number is averaged over 10 different sampled $\theta$'s.
$\mathcal{A}_g$ denotes the Galerkin coarse-grid operator, $\mathcal{A}_{nn}$ is the coarse-grid operator 
obtained from \cref{alg:train}, and SpSA refers to the 
coarse-grid operator from the Sparsified Smooth Aggregation (SpSA) algorithm \cite{treister2015non}.}

\begin{tabular}{c|c|c|c|c  }
 %\hline
 & \multirow{2}{*}{$\xi$} & \multicolumn{3}{c}{$\theta$}  \\
\cline{3-5}
 &  & (${\pi}/{6}$, ${\pi}/{4}$) & (${\pi}/{4}$,${\pi}/{3}$) & 
 (${\pi}/{2}$, ${7\pi}/{12}$) \\
 \hline
 $\mathcal{A}_{g}$ & \multirow{3}{*}{100} & 11.3 & 11.1 & 11.5\\
 ${\mathcal{A}_{nn}}$ & & 16.5 & 16.9 & 19.9 \\
 SpSA & & 40.4 & 36.9 & 11.5 \\
 \hline
 \hline
 $\mathcal{A}_{g}$ & \multirow{3}{*}{200} & 15.9 & 15.3 & 14.5\\
 ${\mathcal{A}_{nn}}$ & & 20.5 & 19.6 & 29.8\\
 SpSA & & 50.4 & 46.7 & 9.1\\
 \hline
 \hline
 $\mathcal{A}_{g}$ & \multirow{3}{*}{300} & 18.1& 21.5 & 17.9 \\
 ${\mathcal{A}_{nn}}$ & & 25.4& 33.1 & 25.6 \\
 SpSA & & 53.0 & 52.1 & 12.3\\
 \hline
  \hline
 $\mathcal{A}_{g}$ & \multirow{3}{*}{400} & 21.1 & 20.2 & 19.9 \\
 ${\mathcal{A}_{nn}}$ & & 27.2 & 30.9 & 26.2 \\
 SpSA & & 59.7 & 56.3 & 12.9\\
 \hline
\end{tabular}
\label{tbl:laplacian-vary-theta-spsa}
\end{table}

\begin{table}[h]
\centering
\caption{The average number of iterations required by 
the GMRES method preconditioned by
3-level AMG methods
with different coarse-grid operators 
for solving \eqref{pde} with 
different sets of PDE parameters. 
The mesh size is $511\times 511$. 
The parameters are selected so that $\theta \in \{\pi/6,\pi/4,\pi/3,5\pi/12\}$ is fixed and $\xi$ is randomly sampled from a uniform distribution in each interval. 
The iteration number is averaged over 10 different sampled $\xi$'s.
$\mathcal{A}_g$ denotes the Galerkin coarse-grid operator, $\mathcal{A}_{nn}$ is the coarse-grid operator 
obtained from \cref{alg:train}, and SpSA refers to the 
coarse-grid operator from the Sparsified Smooth Aggregation (SpSA) algorithm \cite{treister2015non}.}
% \begin{tabular}{c|c|c|c}
% $\theta=\frac{\pi}{4},\xi$ &(0.1,0.13) &(0.13,0.16) &(0.16,0.19)  \\\hline
% $\mathcal{A}_{g}$ &13.1  &11.4 &10.2 \\
% ${\mathcal{A}_{c}}$&17.7 &11.8&10.4 \\

% \end{tabular}

\begin{tabular}{c|c|c|c|c  }
 %\hline
 & \multirow{2}{*}{$\theta$} & \multicolumn{3}{c}{$\xi$}  \\
\cline{3-5}
 &  & ($100, 200$) & ($80, 100$) & 
 ($5, 10$) \\
 \hline
 $\mathcal{A}_{g}$ & \multirow{3}{*}{$\pi / 6$} & 11.6 & 10.4 & 4.4\\
 ${\mathcal{A}_{nn}}$ & & 19.9 & 16.4 & 7.1 \\
 SpSA & &30.3 &27.0  &12.7 \\
 \hline
 \hline
 $\mathcal{A}_{g}$ & \multirow{3}{*}{$\pi/4$} &14.2 &11.3  &4.6\\
 ${\mathcal{A}_{nn}}$ & &18.2 &15.5  &10.0\\
 SpSA & &41.4 &34.7  &13.5\\
 \hline
 \hline
 $\mathcal{A}_{g}$ & \multirow{3}{*}{$\pi/3$} & 11.2 & 10.2 & 4.7 \\
 ${\mathcal{A}_{nn}}$ & & 18.8 & 16.4 & 7.1 \\
 SpSA & &31.7 &28.9  &13.4  \\
 \hline
  \hline
 $\mathcal{A}_{g}$ & \multirow{3}{*}{$5\pi/12$} &11.2 &10.1  &4.8 \\
 ${\mathcal{A}_{nn}}$ & &28.8 &19.1 &9.1 \\
 SpSA & &16.5 &15.1  &9.1\\
 \hline
\end{tabular}
\label{tbl:laplacian-vary-xi-spsa}
\end{table}

\section{Conclusion}
\label{sec-conclusion}

In this work, we propose an ML-based approach 
for computing non-Galerkin coarse-grid operators 
to address the issue of increasing operator complexity in AMG methods by sparsifying the Galerkin operator in different AMG levels. The algorithm consists of two main steps: choosing the sparsity pattern of the stencil and computing the numerical values. 
We employ NNs in both steps and combine their results to construct a  
non-Galerkin coarse-grid
operator with the desired lower complexity.
%\textcolor{red}{RL: summarize a few words on the ML/NN techniques used. Let us wait for 
%Huan to complete Section 3.2}
The NN training algorithm is guided by the AMG convergence theory, ensuring the spectral equivalence of coarse-grid operators with respect to the Galkerkin operator. 
 We showed that spectrally equivalent sparser stencils can be learned by 
advanced ML techniques that exploit
multi-head attention. 

The NN model is trained on parametric PDE problems that cover a wide range of parameters.
The training dataset consists of small-size problems, while the testing problems are significantly larger. Empirical studies conducted on rotated Laplacian and linear elasticity problems provide evidence that the proposed ML method can construct non-Galerkin operators with reduced complexity while maintaining the overall convergence behavior of AMG. A key feature of our method is its ability to generalize to problems of larger sizes and with different PDE parameters that were not encountered in the training. 
This means that the algorithm can effectively handle a wide range of 
problem settings, expanding its practical applicability. By generalizing to new problem instances, the algorithm amortizes the training cost and reduces the need for retraining for every specific problem scenario.

In the future, we plan to extend this work to sparsify unstructured coarse grid operators by exploiting the Graph Convolution Networks (GCNs).  We also plan to explore the  Equivariant Neural Networks \cite{pmlr-v48-cohenc16} to enforce the symmetry in the sparsified coarse-grid operators if the fine level operator is symmetric. Finally, we plan to investigate the real world applications including saddle point system\cite{he2022gdaam}, efficient tensor algebra \cite{9378111, Henderson2018PhenotypingTS,  10.1145/3308558.3313548}, modern generative models \cite{pmlr-v180-cai22a, he2023meddiff, he2021age}, multi-time series analysis techniques \cite{he2023domain, queen2023encoding} to solve time-dependent PDEs.

%We conducted numerical experiments . We found that for both of the challenging problems, our framework can sparsify the coarse grid operator by at least $44\%$ while maintaining the performance of the method. This marks an essential step towards fully boosting the multigrid methods with the help of machine learning models. We also hope our approach can provide new perspectives on improving classical solvers with modern data-driven techniques.

%The multigrid method we are considering in this work is fully defined with matrix-vector products without taking advantage of the underlying geometric information. Therefore it is a legitimate algebraic multigrid method. However, this does not mean that our sparsification algorithm is applicable to every kind of AMG methods. One limitation is that we assume a constant stencil representation for both the fine-level matrix and the coarse-level matrix. Strictly speaking, problems with more complicated meshes would not fit into this framework. Although this does cover a lot of the territory in both practices and theories, it is not in the most generalized version. As such, future work can be devoted to generalizing our framework to unstructured meshes by using graph convolutional neural networks \cite{kipf2016semi}. This will bring up many difficulties for sure. To name a few, one has to consider how complicated it is to set the meshes of the training and testing instances, how to choose the dataset, how to define the feature vector attached to each node, etc. 

\bibliography{reference}
\bibliographystyle{siamplain}

\end{document}